\documentclass[a4paper]{article}

\usepackage{fullpage}

\usepackage{graphics}
\usepackage{epsfig}
\usepackage{amsmath}
\usepackage{amssymb}
\usepackage{bbm}
\usepackage{mathtools}
\usepackage{enumitem}
\usepackage{color} 
\usepackage{float} 

\usepackage[latin1]{inputenc} 
\usepackage[T1]{fontenc} 

\def\Rset{\mathbb{R}}

\def\ds{\displaystyle}
\def\proof{\noindent{\em Proof.} }
\def\qed{\vspace{2mm} \vrule height4pt width3pt depth2pt}

\newtheorem{proposition}{Proposition} 
\newtheorem{assumption}{Assumption} 

\newtheorem{lemma}{Lemma} 
\newtheorem{remark}{Remark}  

\title{Optimal feedback synthesis and minimal time function for the bioremediation of water resources with two patches}
\author{H. Ram\'irez C.$^{1}$, A. Rapaport$^{2}$ and
  V. Riquelme$^{1,2}$\\[2mm]
$^{1}$ Departamento de Ingeniería Matemática \& Centro de Modelamiento
Matemático\\ (UMI 2807, CNRS), Universidad de Chile, Beauchef 851, Casilla 170-3,
Santiago 3, Chile.\\
{\tt hramirez@dim.uchile.cl,vriquelme@dim.uchile.cl}\\[2mm]
$^{2}$ MISTEA, UMR 729 INRA/Supagro, Montpellier, France\\
and MODEMIC, INRA/Inria team, Sophia-Antipolis, France.\\
{\tt rapaport@supagro.inra.fr}}

\date{\today}

\begin{document}

\maketitle

{\em Abstract.}
This paper studies the bioremediation, in minimal time, of a water resource or reservoir using a single continuous bioreactor. The bioreactor is connected to two pumps, at different locations in the reservoir, that pump polluted water and inject back sufficiently clean water with the same flow rate. This leads to a minimal-time optimal control problem where the control variables are related to the inflow rates of both pumps. We obtain a non-convex problem for which it is not possible to directly prove the existence of its solutions. We overcome this difficulty and fully solve the studied problem by applying Pontryagin's principle to the associated generalized control problem. We also obtain explicit bounds on its value function via Hamilton-Jacobi-Bellman techniques.\\


\section{Introduction}

Today, the decontamination of water resources and reservoirs in
natural environments (lakes, lagoons, etc.) and in industrial
frameworks (basin, pools, etc.) is of prime importance. Due to the
availability of drinking water becoming scarce on earth, efforts have
to be made to re-use water and to preserve aquatic resources. To this
end, biological treatment is a convenient way to extract organic or
soluble matter from water. The basic principle is to use biotic agents
(generally micro-organisms) that convert the pollutant until the
concentration in the reservoir decreases to an acceptable
level. Typically, the treatment is performed with the help of
continuously stirred or fed-batch bioreactors. Numerous studies have
been devoted to this subject over the past 40 years (see, for
instance, 
\cite{AGK72,AKG71,BEMA03,GHR08,KH00,M99,SCNBV04,SV02,SKLB08,SPB03,VB98}).
 
The following main types of procedure are usually considered:
\begin{itemize}
\item The direct introduction of the biotic agents to the reservoir. This solution could lead to the eutrophication of the resource.
\item  The draining of the reservoir to a dedicated bioreactor and the
filling back of the water after treatment. This solution attempts to eradicate various forms of life supported by the water resource, that cannot survive without water (such as fish, algae, etc.).
\end{itemize}
Alternatively, one can consider a side bioreactor that continuously treats the water pumped from the reservoir and that injects it back with the same flow rate so that the volume of the reservoir remains constant at all time. At the output of the bioreactor, a settler separates biomass from the water so that no biomass is introduced in the resource. Such an operating procedure is typically used for water purification of culture basins in aquaculture \cite{CADBV07,EKVHK06,P03}.

The choice of the flow rate presents a trade-off between the speed at which the water is treated and the quality of decontamination. Recently, minimal-time control problems with simple spatial representations have been formulated and addressed \cite{GHRR11}. Under the assumption that the resource is perfectly mixed, an optimal state-feedback that depends on the characteristics of the micro-organisms and on the on-line measurement of the pollutant concentration has been derived. Later, an extension with a more realistic spatial representation was proposed in \cite{GRRR12} that considers two perfectly-mixed zones in the resource: an ``active'' zone, where the treatment of the pollutant is the most effective, and a more confined or ``dead'' zone that communicates with the active zone by diffusion of the pollutant.
It has been shown that the optimal feedback obtained for the
perfectly mixed case is also
optimal when one applies it on the pollutant concentration in the active
zone only. The fact that this controller does not require knowledge of
the size of the dead zone or of the value of the diffusion parameter,
neither of the online measurement of the pollutant in the dead zone,
is a remarkable property. Nevertheless, the minimal time is impacted
by the characteristics of the confinement. 

\begin{figure}[h]
\begin{center}
\includegraphics[scale=0.35]{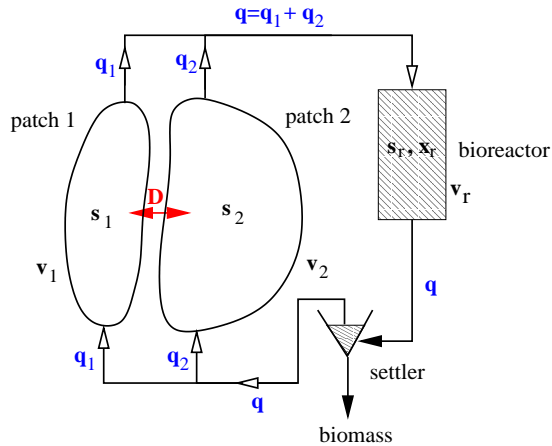}
\caption{Modeling scheme of the treatment of two interconnected
  patches (definitions of the variables and parameters are given
  in Section \ref{sec_def})}\label{fig1}
\end{center}
\end{figure}
In the present work, we consider that the
  treatment of the water resource can be split into two zones i.e.
  the water is extracted from the resource at
  two different points (instead of one), and the treated water
  returns to the resource (with the same flows) at two different
  locations. A diffusion makes connection between the zones (see
  Fig. \ref{fig1}).
Such a division into two patches can represent real situations such as:
\begin{itemize}
\item natural environments where water tables or lagoons are
  connected together by a small communication path (this modeling
  covers also the particular case of a null diffusion when one has to treat
  two independent volumes),
\item resource hydrodynamics that reveal influence zones for each
  pumping devices, depending on the locations of the extraction and
  return points,
\item accidental pollution as an homogeneous strain diffusing into the
  complementary part of the resource.
\end{itemize}
The control problem consists in choosing dynamically the total flow
rate $q$ and the flow distribution $q_{1}$, $q_{2}$ between the two patches, with the objective of
having both of them decontaminated in minimal time.
Notice that a particular strategy consists in having all the time a
flow distribution entirely with one zone, which amounts to the former problem
with active and dead zones mentioned above. We study here the benefit
of switching dynamically the treatment to the other patch or treating
simultaneously both patches.
The associated minimal-time problem is significantly more complex, because there are two controls and the velocity set is non-convex (this is shown in the next Section). Indeed, it is necessary to use different techniques to address the cases of non-null diffusion between the two zones and the limiting case of null diffusion between the two zones.

The paper is organized as follows. In the next section, definitions
and assumptions are presented. In Section \ref{sec_relaxed}, properties
of the optimization problem with relaxed controls and non-null
diffusion are investigated. In Section \ref{sec_feedback}, the optimal
control strategy for the original problem with non-null diffusion is
given and proven. In Section \ref{sec_valuefunction}, we address the
particular case of null diffusion and we provide explicit bounds on
the minimal-time function. Finally, we show numerical computations
that illustrate the theoretical results, and give concluding remarks.

\section{Definitions and preliminaries}
\label{sec_def}

In what follows, we denote by $\Rset$ the set of real
    numbers, $\Rset_{+}$ and $\Rset_{+}^{\star}$ the sets of
    non-negative and positive real numbers respectively.
    Analogously, $\Rset_{-}$ and $\Rset_{-}^{\star}$ are the sets of
    non-positive and negative real numbers respectively. We set also
    $\Rset_{+}^{2}=\Rset_{+}\times\Rset_{+}$ and
    $\Rset_{-}^{2}=\Rset_{-}\times\Rset_{-}$.

The time evolution of the concentrations $s_{i}$ ($i=1,2$) of pollutants in the two patches are given by the equations
\begin{equation}
\label{dyn-resource}
\left\{\begin{array}{lll}
\ds \frac{ds_{1}}{dt} & = & \ds
  \frac{q_{1}}{v_{1}}(s_{{\rm{r}}}-s_{1})+\frac{D}{v_{1}}(s_{2}-s_{1}) \ ,\\[3mm]
\ds \frac{ds_{2}}{dt} & = & \ds
\frac{q_{2}}{v_{2}}(s_{{\rm{r}}}-s_{2})+\frac{D}{v_{2}}(s_{1}-s_{2}) \ ,
\end{array}\right.
\end{equation}
where the volumes $v_{i}$ ($i=1,2$) are assumed to be constant and $D$ denotes the diffusion coefficient of the pollutant between the two zones. The control variables are the flow rates $q_{i}$ of the pumps in each zone, which bring water with a low pollutant concentration $s_{{\rm{r}}}$ from the bioreactor and remove water with a pollutant concentration $s_{i}$ from each zone $i$, with the same flow rates $q_{i}$.

The concentration $s_{{\rm{r}}}$ at the output of the bioreactor is linked to the total flow rate $q=q_{1}+q_{2}$ by the usual chemostat model:
\begin{equation}
\label{dyn-bioreactor}
\left\{\begin{array}{lll}
\ds \frac{ds_{{\rm{r}}}}{dt} & = &
\ds -\mu(s_{{\rm{r}}})x_{{\rm{r}}}+\frac{q}{v_{{\rm{r}}}}(s_{in}-s_{{\rm{r}}}) \ ,\\[3mm]
\ds \frac{dx_{{\rm{r}}}}{dt}  & =  & \ds \mu(s_{{\rm{r}}})x_{{\rm{r}}}-\frac{q}{v_{{\rm{r}}}}x_{{\rm{r}}} \ ,
\end{array}\right.
\end{equation}
where $x_{{\rm{r}}}$ is the biomass concentration, $v_{{\rm{r}}}$ is the volume of the bioreactor and $\mu(\cdot)$ is the
specific growth rate of the bacteria (without a loss of generality we
assume that the yield coefficient is equal to one). 
These equations describe the dynamics of a bacterial growth
  consuming a substrate that is constantly fed in a tank of constant
  volume (see for instance \cite{SW95}).
The input concentration $s_{in}$ is given here by the combination of the
concentrations of the water extracted from the two zones:
\begin{equation}
\label{coupling}
s_{in}=\frac{q_{1}s_{1}+q_{2}s_{2}}{q_{1}+q_{2}} \ .
\end{equation}
We assume that the output of the bioreactor is filtered by a
settler, that we assume to be perfect, so that the water that returns to the resource is biomass
free (see \cite{DF13a,DF13b} for considerations of settler modeling and
  conditions that ensure the stability of the desired steady-state of the settler).\\

The target to be reached in the minimal time is defined by a threshold $\underline s>0$ of the pollutant concentrations, that is
\begin{equation}
\label{s_in}
{\cal T}=\left\{ s=(s_1,s_2) \in\Rset_{+}^{2} \, \vert \,
\max(s_{1},s_{2})\leq\underline s \right\} \ .
\end{equation}
In the paper, we shall denote $t_{f}$ as the first time at which a
trajectory reaches the target (when it exists).\\

We make the usual assumptions on the growth function $\mu(\cdot)$
in absence of inhibition.

\begin{assumption}
\label{H0}
$\mu(\cdot)$ is a
$C^{1}$ increasing concave function defined on $\Rset_{+}$ with $\mu(0)=0$.
\end{assumption}

Under this last assumption, we recall that under a constant $s_{in}$, the dynamics
(\ref{dyn-bioreactor}) admit a unique positive equilibrium $(s_{{\rm{r}}}^{\star},x_{{\rm{r}}}^{\star})$ that is globally asymptotically stable on the domain $\Rset_{+}\times\Rset_{+}^{\star}$ provided that the condition $q/v_{{\rm{r}}}\leq \mu(s_{in})$ is satisfied (see, for instance, \cite{SW95}). Then, $s_{{\rm{r}}}^{\star}$ is defined as the unique
solution of $\mu(s_{{\rm{r}}}^{\star})=q/v_{{\rm{r}}}$ and $x_{{\rm{r}}}^{\star}=s_{in}-s_{{\rm{r}}}^{\star}$.
Consequently, considering expression (\ref{coupling}), the controls $q_{1}$ and $q_{2}$ are chosen such that
\begin{equation}
q_{1}+q_{2} \leq
v_{{\rm{r}}}\mu\left(\frac{q_{1}s_{1}+q_{2}s_{2}}{q_{1}+q_{2}}\right) \ .
\end{equation}

We assume that the resource to be treated is very large. This amounts to considering that the bioreactor is small compared to both zones of the resource.

\begin{assumption}
\label{H1}
$v_{1}$ and $v_{2}$ are large compared to $v_{{\rm{r}}}$.
\end{assumption}

Let us define $\alpha=q_{1}/q$, $r=v_{1}/(v_1+v_2)$, $d=D/v_{{\rm{r}}}$, and $\epsilon=v_{{\rm{r}}}/(v_1+v_2)$. Then, the coupled dynamics
(\ref{dyn-resource})-(\ref{dyn-bioreactor}) with (\ref{coupling}) can
be written in the slow-fast form
\begin{equation}
\label{slowfast}
\left\{\begin{array}{lll}
\ds \frac{ds_{{\rm{r}}}}{dt} & = & \ds
-\mu(s_{{\rm{r}}})x_{{\rm{r}}}+\frac{q}{v_{{\rm{r}}}}(\alpha s_{1}+(1-\alpha)s_{2}-s_{{\rm{r}}}) \
, \\[3mm]
\ds \frac{dx_{{\rm{r}}}}{dt} & = &  \ds \mu(s_{{\rm{r}}})x_{{\rm{r}}}-\frac{q}{v_{{\rm{r}}}}x_{{\rm{r}}} \
, \\[3mm]
\ds  \frac{ds_{1}}{dt} & = &  \ds \epsilon\left(\frac{\alpha}{r}\frac{q}{v_{{\rm{r}}}}(s_{{\rm{r}}}-s_{1})+\frac{d}{r}(s_{2}-s_{1})\right) \
, \\[3mm]
\ds  \frac{ds_{2}}{dt} & = & \ds
\epsilon\left(\frac{1-\alpha}{1-r}\frac{q}{v_{{\rm{r}}}}(s_{{\rm{r}}}-s_{2})+\frac{d}{1-r}(s_{1}-s_{2})\right) \
.
\end{array}\right.
\end{equation}
Provided that the initial conditions of the variables $(s_{{\rm{r}}},x_{{\rm{r}}})$
belong to $\Rset_{+}\times\Rset_{+}^{\star}$, applying Tikonov's Theorem (see
for instance \cite{K01}), the dynamics of the slow variables $(s_{1},s_{2})$ can be approached using the reduced dynamics
\begin{equation}
\label{reduced}
\left\{\begin{array}{lllll}
\ds \dot s_{1} & = & \ds \frac{ds_{1}}{d\tau} & = & \ds
  \frac{\alpha}{r}\mu(s_{{\rm{r}}}^{\star})(s_{{\rm{r}}}^{\star}-s_{1})+\frac{d}{r}(s_{2}-s_{1})
  \ , \\[3mm]
\ds \dot s_{2} & = & \ds \frac{ds_{2}}{d\tau} & = & \ds 
\frac{1-\alpha}{1-r}\mu(s_{{\rm{r}}}^{\star})(s_{{\rm{r}}}^{\star}-s_{2})+\frac{d}{1-r}(s_{1}-s_{2})
\ 
\end{array}\right.
\end{equation}
in the time scale $\tau=\epsilon t$. In this formulation, the quasi-steady-state
concentration $s_{{\rm{r}}}^{\star}$ of the bioreactor can be considered as a
control variable that takes values in $[0,\alpha
s_{1}+(1-\alpha)s_{2}]$, which is equivalent to choosing $q \in
[0,v_{{\rm{r}}}\mu(\alpha s_{1}+(1-\alpha)s_{2})]$ when Assumption \ref{H0}
is satisfied. 
In the following, we shall consider the optimal control for the
  reduced dynamics only. Nevertheless, we give some
  properties of the optimal feedback for the reduced dynamics when
  applied to the un-reduced one, in Section \ref{sec_feedback} (Remark
  \ref{rem_unreduced}) and
  Appendix. 

Notice that the control problem can be reformulated with the controls $u=(\alpha, s_{{\rm{r}}}^{\star})$ that belong to the state-dependent control set
\begin{equation}
U(s)=\left\{ (\alpha,s_{{\rm{r}}}^{\star}) \, \vert \, \alpha \in [0,1], \;
s_{{\rm{r}}}^{\star}\in[0,\alpha s_{1}+(1-\alpha)s_{2}]\right\}
\end{equation}
equivalently to controls $q_{1}$ and $q_{2}$. In what follows, a measurable function $u(\cdot)$ such that $u(t)\in U(s(t))$ for all $t$ is called an \emph{admissible control}.

\begin{lemma}
\label{lemma1}
The domain $\Rset_{+}^{2}$ is positively invariant by the dynamics (\ref{reduced}) for any admissible controls $u(\cdot)$, and any trajectory is bounded. Furthermore, the target ${\cal T}$ is reachable in a finite time from any initial condition in $\Rset_{+}^{2}$.
\end{lemma}

\proof  
For $s_{1}=0$ and $s_{2}\geq 0$, one has $\dot s_{1}\geq 0$. Similarly, one has $\dot s_{2}\geq 0$ when $s_{1}\geq 0$ and $s_{2}=0$. By the uniqueness of the solutions of (\ref{reduced}) for measurable controls $u(\cdot)$, we deduce that $\Rset_{+}^{2}$ is invariant. From equations (\ref{reduced}), one can write
\[
r\dot s_{1}+(1-r)\dot s_{2}=\mu(s_{{\rm{r}}}^{\star})(s_{{\rm{r}}}^{\star}-(\alpha
s_{1}+(1-\alpha)s_{2}))\leq 0
\]
for any admissible controls. One then deduces
\[
rs_{1}(t)+(1-r)s_{2}(t)\leq M_{0}=rs_{1}(0)+(1-r)s_{2}(0), \quad \forall t
\geq 0,
\]
which provides the boundedness of the trajectories.\\

Consider the feedback strategy
\[
\alpha=r \ , \; s_{{\rm{r}}}^{\star}=\frac{rs_{1}+(1-r)s_{2}}{2},
\]
and we write the dynamics of $m=rs_{1}+(1-r)s_{2}$ as follows:
\[
\dot m=-\mu\left(\frac{m}{2}\right)\frac{m}{2} < 0 \ , \forall m>0 \ .
\]
Then, from any initial condition in $\Rset_{+}^{2}$, the solution
$m(t)$ tends to $0$ when $t$ tends to infinity. Therefore, $m(\cdot)$ reaches the set $[0,\min(r,1-r)\underline s]$ in a finite time, which guarantees that
$s=(s_1,s_2)$ belongs to ${\cal T}$ at that time.\qed

\bigskip

For simplicity, we define the function
\begin{equation}
\beta(\sigma,s_{{\rm{r}}}^{\star})=\mu(s_{{\rm{r}}}^{\star})(\sigma-s_{{\rm{r}}}^{\star}) 
\end{equation}
so that the dynamics (\ref{reduced}) can be written in the more compact form
\begin{equation}
\label{dynamics}
\dot s = F(s,u)+dG(s)
\end{equation}
where $F(\cdot)$ and $G(\cdot)$ are defined as follows:
\[
F(s,(\alpha,s_{{\rm{r}}}^{\star}))=-\left[\begin{array}{c}
\ds \frac{\alpha}{r}\beta(s_{1},s_{{\rm{r}}}^{\star})\\[2mm]
\ds \frac{1-\alpha}{1-r}\beta(s_{2},s_{{\rm{r}}}^{\star})
\end{array}\right]
, \quad
G(s)=\left[\begin{array}{c}
\ds \frac{s_{2}-s_{1}}{r}\\[2mm]\ds \frac{s_{1}-s_{2}}{1-r}
\end{array}\right] \ .
\]
The dynamics can be equivalently expressed in terms of controls $v=(\alpha,\zeta)$ that belong to the state-independent set $V=[0,1]^{2}$ with
the dynamics
\begin{equation}
\label{dynamics_with_u}
\dot s = F(s,(\alpha,\zeta(\alpha s_{1}+(1-\alpha)s_{2})))+dG(s)
\end{equation}
which satisfy the usual regularity conditions for applying Pontryagin's Maximum Principle for deriving necessary optimality conditions. One can notice that the velocity set of the dynamics (\ref{dynamics_with_u}) is not everywhere convex. 
Consequently, one cannot guarantee a priori the existence of an optimal control $v(\cdot)$ in the set of time-measurable functions that take values in $V$ but that are among relaxed controls (see, for instance, \cite[Sec. 2.7]{V00}). For convenience, we shall keep the formulation of the problem with controls $u$.
Because for any $s$ the sets $\cup_{u \in U(s)}F(s,u)$ are two-dimensional connected,
the corresponding convexified dynamics can be written as follows (see \cite[Th. 2.29]{RW98}):
\begin{equation}
\label{dynamics_convexified}
\dot s = \tilde F(s,\tilde u)+dG(s)
\end{equation}
with
\begin{equation}
\label{F_convexified}
\tilde F(s,\tilde u)=p F(s,u_{a})+(1-p)F(s,u_{b})
\end{equation}
where the relaxed controls $\tilde u=(u_{a},u_{b},p)=(\alpha_{a},s_{{\rm{r}}a}^{\star},\alpha_{b},s_{{\rm{r}}b}^{\star},p)$ belong to the set 
\[
\tilde U(s)=U(s)^{2}\times[0,1].
\] 

In the next section, we show that the relaxed problem admits an optimal solution that is also a solution of the original (non-relaxed) problem.

\section{Study of the relaxed problem}
\label{sec_relaxed}

Throughout this section, we assume that the parameter $d$ is positive. The particular case of $d=0$ will be considered later in Section \ref{sec_valuefunction}.
Let us write the Hamiltonian of the relaxed problem
\begin{equation}
\label{Hamiltonian_relaxed}
\begin{multlined}
\tilde H(s,\lambda,(\alpha_{a},s_{{\rm{r}}a}^{\star},\alpha_{b},s_{{\rm{r}}b}^{\star},p))=\\  
-1+p\,Q(s,\lambda,(\alpha_{a},s_{{\rm{r}}a}^{\star}))+(1-p)\,Q(s,\lambda,(\alpha_{b},s_{{\rm{r}}b}^{\star}))+d(s_{2}-s_{1})\left(\frac{\lambda_{1}}{r}-\frac{\lambda_{2}}{1-r}\right)
\end{multlined}
\end{equation}
which is to be maximized w.r.t. $(\alpha_{a},s_{{\rm{r}}a}^{\star},\alpha_{b},s_{{\rm{r}}b}^{\star},p)\in \tilde U(s)$, where $\lambda=(\lambda_1,\lambda_2)$, and we have defined, for convenience, the function
\begin{equation}
\label{defQ}
Q(s,\lambda,(\alpha,s_{{\rm{r}}}^{\star}))=
-\left(\alpha\frac{\lambda_{1}}{r}\beta(s_{1},s_{{\rm{r}}}^{\star})+(1-\alpha)\frac{\lambda_{2}}{1-r}\beta(s_{2},s_{{\rm{r}}}^{\star})\right)
\ .
\end{equation}
The adjoint equations are
\begin{equation}
\label{dyn-adjoint}
\left\{\begin{array}{lll}
\ds \dot\lambda_{1} & = & \ds \lambda_{1}\left(p\frac{\alpha_{a}}{r}\mu(s_{{\rm{r}}a}^{\star})+(1-p)\frac{\alpha_{b}}{r}\mu(s_{{\rm{r}}b}^{\star})
+\frac{d}{r}\right)-\lambda_{2}\frac{d}{1-r} \ ,\\[3mm]
\ds \dot\lambda_{2} & = & \ds -\lambda_{1}\frac{d}{r}+
\lambda_{2}\left(p\frac{1-\alpha_{a}}{1-r}\mu(s_{{\rm{r}}a}^{\star})
+(1-p)\frac{1-\alpha_{b}}{1-r}\mu(s_{{\rm{r}}b}^{\star})
+\frac{d}{1-r}\right)
\ ,
\end{array}\right.
\end{equation}
with the following transversality conditions
\begin{equation}
\label{transversality}
\left\{\begin{array}{lll}
s_{1}(t_{f})<\underline{s}, \; s_{2}(t_{f})=\underline{s}  &
  \Rightarrow & \lambda_{1}(t_{f})=0, \; \lambda_{2}(t_{f})<0 \ ,\\
s_{1}(t_{f})=\underline{s}, \; s_{2}(t_{f})<\underline{s} & \Rightarrow
& \lambda_{1}(t_{f})<0, \; \lambda_{2}(t_{f})=0 \ , \\
s_{1}(t_{f})=\underline{s}, \; s_{2}(t_{f})=\underline{s} & \Rightarrow
& \lambda_{1}(t_{f})\leq0, \; \lambda_{2}(t_{f})\leq 0 \mbox{ with }
 \lambda(t_{f})\neq 0 \ .
\end{array}\right.
\end{equation}

As usual, a triple $(s(\cdot),\lambda(\cdot),\tilde u^{\star}(\cdot))$ satisfying \eqref{dynamics_convexified}, \eqref{dyn-adjoint}, \eqref{transversality}, and 
\begin{equation}\label{eq:maxRelHam}
\tilde H(s(t),\lambda(t),\tilde u^{\star}(t)) = \max_{\tilde u\in\tilde U(s(t))} \tilde H(s(t),\lambda(t),\tilde u)
\end{equation}
is called an \emph{admissible extremal}.

\begin{lemma}
\label{lemma2}
Along any admissible extremal, one has $\lambda_{i}(t)<0$ ($i=1,2$) for any $t<t_{f}$.
\end{lemma}

\proof
If one writes the adjoint equations \eqref{dyn-adjoint} as
$\dot \lambda_{i}=\phi_{i}(t,\lambda_{1},\lambda_{2})$ ($i=1,2$), one
can notice that the partial derivatives $\partial_{j}\phi_{i}$ ($i\neq
j$) are non-positive. From the theory of monotone dynamical systems
(see for instance \cite{S95}), the dynamics \eqref{dyn-adjoint} is thus competitive
or, equivalently, cooperative in backward time.
As the transversality conditions (\ref{transversality}) gives
$\lambda_{i}(t_{f})\leq 0$ ($i=1,2$), we deduce by the property
of monotone dynamics that one should have $\lambda_{i}(t)\leq 0$
($i=1,2$) for any $t\leq t_{f}$. Moreover, $\lambda=0$ is an
equilibrium of  \eqref{dyn-adjoint} and $\lambda(t_{f})$ has to be
different from $0$ at any time $t\leq t_{f}$. Then, $\lambda_{i}(t)$
($i=1,2$) cannot be simultaneously equal to zero.
If there exists $t<t_{f}$ and $i \in \{1, 2\}$ such that
$\lambda_i(t)=0$, then one should have $\lambda_j(t)<0$ for $j\neq
i$. However, $d>0$ implies $\dot \lambda_i(t)>0$, thus obtaining a contradiction with $\lambda_i\leq 0$ for any time.
\qed

\bigskip

For the following, we consider the function 
\begin{equation}
\label{gamma}
\gamma(\sigma)=\max_{s_{{\rm{r}}}^{\star}\geq 0}\beta(\sigma,s_{{\rm{r}}}^{\star}), \quad
\sigma > 0 \ ,
\end{equation}
which satisfies the following property:

\begin{lemma} 
\label{lemma3}
Under Assumption \ref{H0}, for any $\sigma > 0$, there exists a unique $\hat s_{{\rm{r}}}^{\star}(\sigma) \in (0,\sigma)$ that realizes the maximum in (\ref{gamma}). Furthermore, the function $\gamma(\cdot)$ is differentiable and increasing with
\begin{equation}
\label{gamma_prime}
\gamma^{\prime}(\sigma)=\mu(\hat s_{{\rm{r}}}^{\star}(\sigma)) \ .
\end{equation}
\end{lemma}

\proof
Consider the function $\varphi: (\sigma,w)\in \Rset_{+}\times [0,1] \mapsto
\beta(\sigma,w\sigma)$ and the partial function $\varphi_{\sigma}: w \in
[0,1] \mapsto \varphi(\sigma,w)$ for fixed $\sigma>0$. Notice that $\varphi_{\sigma}(0)=\varphi_{\sigma}(1)=0$ and that $\varphi_{\sigma}(w)>0$ for $w \in (0,1)$. 
Simple calculation gives $\varphi_{\sigma}^{\prime\prime}(w)=\mu^{\prime\prime}(w\sigma)(1-w)\sigma^{3}-2\mu^{\prime}(w\sigma)\sigma^{2}$,
which is negative. Therefore, $\varphi_{\sigma}(\cdot)$ is a strictly concave
function on $[0,1]$ and consequently admits a unique maximum
$w^{\star}_{\sigma}$ on $[0,1]$. We conclude that $w^{\star}_{\sigma}$
belongs to $(0,1)$ or, equivalently, that the maximum of $s_{{\rm{r}}}^{\star} \mapsto \beta(\sigma,s_{{\rm{r}}}^{\star})$ is realized for a unique  $\hat s_{{\rm{r}}}^{\star}(\sigma)=w^{\star}_{\sigma}\sigma$ in $(0,\sigma)$.

Furthermore, one has $\varphi_{\sigma}^{\prime}(w)=\sigma\mu^{\prime}(w\sigma)(\sigma-w\sigma)-\sigma\mu(w\sigma)$, and the necessary optimality condition $\varphi_{\sigma}^{\prime}(w^{\star}_{\sigma})=0$ gives the equality
\begin{equation}
\label{property}
\mu(\hat s_{{\rm{r}}}^{\star}(\sigma))=\mu^{\prime}(\hat s_{{\rm{r}}}^{\star}(\sigma))(\sigma-\hat
s_{{\rm{r}}}^{\star}(\sigma)) \ .
\end{equation}
Simple calculation shows that for each $w\in[0,1]$, the function $\sigma\mapsto\varphi(\sigma,w)$ is convex. Because the maximizer $w_{\sigma}^{\star}$ of $\varphi_{\sigma}(\cdot)$ is unique for any $\sigma$, one can apply the rules of differentiability of pointwise maxima (see, for instance, \cite[Chap. 2.8]{C87}), which state that the function $\gamma(\sigma)=\max_{w \in   [0,1]}\varphi(\sigma,w)$ is differentiable with
\[
\gamma^{\prime}(\sigma)=\frac{\partial \varphi}{\partial
  \sigma}(\sigma,w^{\star}_{\sigma}\sigma)
=w_{\sigma}^{\star}\mu^{\prime}(w_{\sigma}^{\star}\sigma)(\sigma-w_{\sigma}^{\star}\sigma)+\mu(w_{\sigma}^{\star}\sigma)(1-w_{\sigma}^{\star})
\ .
\]
Equation (\ref{property}) provides the simpler expression (\ref{gamma_prime}), which shows that $\gamma(\cdot)$ is increasing. 
\qed

\bigskip

We now consider the variable
\begin{equation}
\label{def_eta}
\eta=\frac{-\lambda_{1}}{r}\,\gamma(s_{1})-\frac{-\lambda_{2}}{1-r}\,\gamma(s_{2})
\end{equation}
which will play the role of a {\em switching function}. Notice that this is
not the usual switching function of problems with linear dynamics w.r.t. a scalar control because our problem has two controls $\alpha$ and $s_{{\rm{r}}}^{\star}$ that cannot be separated, and the second control acts non-linearly in the dynamics.

\begin{lemma} 
\label{lemma4}
For fixed $(s,\lambda)\in\Rset_{+}^{2}\times\Rset_{-}^{2}$, the pairs $u^{\star}=(\alpha,s_r^{\star}) \in U(s)$ that maximize the function $Q(s,\lambda,\cdot)$ are the following:
\begin{enumerate}[label=\roman{*}.]
\item when $\eta>0$: $u^{\star}=(1,\hat s_{{\rm{r}}}^{\star}(s_{1}))$,
\item when $\eta<0$: $u^{\star}=(0,\hat s_{{\rm{r}}}^{\star}(s_{2}))$,
\item when $\eta=0$ and $s_{1}= s_{2}$: 
$u^{\star}\in [0,1]\times \{\hat s_{{\rm{r}}}^{\star}\}$
where $\hat s_{{\rm{r}}}^{\star}=\hat s_{{\rm{r}}}^{\star}(s_{1})=\hat s_{{\rm{r}}}^{\star}(s_{2})$,
\item  when $\eta=0$ and $s_{1}\neq s_{2}$:
  $u^{\star}=(1,\hat s_{{\rm{r}}}^{\star}(s_{1}))$ or  $u^{\star}=(0,\hat s_{{\rm{r}}}^{\star}(s_{2}))$.
\end{enumerate}
\end{lemma}

\proof
When $\eta>0$, one can write, using Lemma \ref{lemma3} and $\lambda_{1},\lambda_2<0$,
\[
\begin{array}{lll}
Q(s,\lambda,(1,\hat s_{{\rm{r}}}^{\star}(s_{1}))) & = &
\ds \frac{-\lambda_{1}}{r}\gamma(s_{1})\\
 & > &
 \ds \alpha\frac{-\lambda_{1}}{r}\gamma(s_{1})+(1-\alpha)\frac{-\lambda_{2}}{1-r}\gamma(s_{2})
 \ ,
 \; \forall \alpha \in [0,1)\\
 & \geq &
 \ds \alpha\frac{-\lambda_{1}}{r}\beta(s_{1},s_{{\rm{r}}}^{\star})+(1-\alpha)\frac{-\lambda_{2}}{1-r}\beta(s_{2},s_{{\rm{r}}}^{\star})
 \ ,
 \; \forall \alpha \in [0,1), \; \forall s_{{\rm{r}}}^{\star} \in [0,\alpha
 s_{1}+(1-\alpha)s_{2}],
\\
 & \geq &
\ds Q(s,\lambda,(\alpha, s_{{\rm{r}}}^{\star}))
 \ ,
 \; \forall \alpha \in [0,1), \; \forall s_{{\rm{r}}}^{\star} \in [0,\alpha
 s_{1}+(1-\alpha)s_{2}],
\end{array}
\]
and for $\alpha=1$, one has $Q(s,\lambda,(1,\hat s_{{\rm{r}}}^{\star}(s_{1}))) > Q(s,\lambda,(1,s_{{\rm{r}}}^{\star})) \ , \forall s_{{\rm{r}}}^{\star} \neq \hat s_{{\rm{r}}}^{\star}(s_{1})$. Therefore, the maximum of $Q(s,\lambda,\cdot)$ is reached for the unique pair $(\alpha,s_{{\rm{r}}}^{\star})=(1,s_{{\rm{r}}}^{\star}(s_{1}))$.\\

Similarly, when $\eta<0$, one can show that the unique maximum is $(\alpha,s_{{\rm{r}}}^{\star})=(0,s_{{\rm{r}}}^{\star}(s_{2}))$.\\

When $\eta=0$, one has
\[
\frac{-\lambda_{1}}{r}\gamma(s_{1})=\frac{-\lambda_{2}}{1-r}\gamma(s_{2})
> Q(s,\lambda,(\alpha,s_{{\rm{r}}}^{\star}))
\ , \forall \alpha \in [0,1] \ , \forall s_{{\rm{r}}}^{\star}\notin \{ \hat
s_{{\rm{r}}}^{\star}(s_{1}), \hat s_{{\rm{r}}}^{\star}(s_{2}) \}.
\]
If $s_{1}=s_{2}$, one necessarily has $\lambda_{1}/r=\lambda_{2}/(1-r)\neq 0$, and thus,
\[
Q(s,\lambda,(\alpha,s_{{\rm{r}}}^{\star}))=\frac{-\lambda_{1}}{r}\beta(s_{1},s_{{\rm{r}}}^{\star})<
\frac{-\lambda_{1}}{r}\gamma(s_{1})=Q(s,\lambda,(\alpha,\hat
s_{{\rm{r}}}^{\star}(s_{1}))) \ , \forall s_{1}^{\star} \neq \hat
s_{{\rm{r}}}^{\star}(s_{1}),
\]
for any $\alpha \in [0,1]$. The optimal $s_{{\rm{r}}}^{\star}$ is necessarily equal to $\hat s_{{\rm{r}}}^{\star}(s_{1})=\hat s_{{\rm{r}}}^{\star}(s_{2})$.\\

If $s_{1}\neq s_{2}$, one has $\tilde s_{{\rm{r}}}^{\star}(s_{1}) \neq \hat s_{{\rm{r}}}^{\star}(s_{2})$, and
consequently, using Lemma \ref{lemma3} and the fact that $\lambda_{1}$ and $\lambda_{2}$ are both negative,
\[
Q(s,\lambda,(\alpha,\hat s_{{\rm{r}}}^{\star}(s_{1}))) = 
\alpha\frac{-\lambda_{1}}{r}\gamma(s_{1})+(1-\alpha)\frac{-\lambda_{2}}{1-r}\beta(s_{2},\hat
s_{{\rm{r}}}^{\star}(s_{1}))
< \frac{-\lambda_{1}}{r}\gamma(s_{1}) \ , \forall \alpha \in[0,1)
\]
\[
Q(s,\lambda,(\alpha,\hat s_{{\rm{r}}}^{\star}(s_{2}))) = 
\alpha\frac{-\lambda_{1}}{r}\beta(s_{1},\hat
s_{{\rm{r}}}^{\star}(s_{2}))+(1-\alpha)\frac{-\lambda_{2}}{1-r}\gamma(s_{2})
< \frac{-\lambda_{2}}{1-r}\gamma(s_{2}) \ , \forall \alpha \in(0,1]
\]
Then,  $(\alpha,s_{{\rm{r}}}^{\star})=(1,\hat s_{{\rm{r}}}^{\star}(s_{1}))$ and $(\alpha,s_{{\rm{r}}}^{\star})=(0,\hat s_{{\rm{r}}}^{\star}(s_{2}))$ are the only two pairs that maximize $Q(s,\lambda,\cdot)$.
\qed

\begin{proposition}
\label{prop_tilde_u}
At almost any time, an optimal control $\tilde u^{\star}$ of the relaxed problem satisfies the following property:
\begin{enumerate}
\item when $\eta \neq 0$ or $s_{1}=s_{2}$, one has
$\tilde F(s,\tilde u^{\star})=F(s,u^{\star})$, where $u^{\star}$ is given by Lemma \ref{lemma4} i.-ii.-iii.

\item when $\eta = 0$ and $s_{1}\neq s_{2}$, one has
\begin{equation}
\label{chattering}
\begin{multlined}
\tilde u^{\star} \in \{(1,\hat s_{{\rm{r}}}^{\star}(s_{1})),(0,\hat
s_{{\rm{r}}}^{\star}(s_{2}))\}\times U(s)\times\{1\}
\; \cup \;
U(s)\times\{(1,\hat s_{{\rm{r}}}^{\star}(s_{1})),(0,\hat
s_{{\rm{r}}}^{\star}(s_{2}))\}\times\{0\}\\
\; \cup \;
\{(1,\hat s_{{\rm{r}}}^{\star}(s_{1}),0,\hat s_{{\rm{r}}}^{\star}(s_{2}))\}\times[0,1]
\; \cup \;
\{(0,\hat s_{{\rm{r}}}^{\star}(s_{2}),1,\hat s_{{\rm{r}}}^{\star}(s_{1}))\}\times[0,1].
\end{multlined} 
\end{equation}
\end{enumerate}
\end{proposition}

\proof
According to Pontryagin's Maximum Principle, an optimal control $\tilde u=(u_{a},u_{b},p)$ has to maximize for a.e. time the Hamiltonian $\tilde H$ given in (\ref{Hamiltonian_relaxed}) or, equivalently, the quantity
\[
(u_{a},u_{b},p) \longmapsto \tilde Q(s,\lambda,(u_{a},u_{b},p))=pQ(s,\lambda,u_{a})+(1-p)Q(s,\lambda,u_{b})
\]
where $\lambda_{1}$ and $\lambda_{2}$ are negative (from Lemma \ref{lemma2}).
Let us consider the maximization of the function $Q(s,\lambda,\cdot)$ characterized by Lemma \ref{lemma4}.

In cases i and ii, the function $Q(s,\lambda,\cdot)$ admits a unique maximizer $u^{\star}$. Thus, $\tilde Q(s,\lambda,\cdot)$ is maximized for $u_{a}=u^{\star}$ with $p=1$ independent of $u_{b}$ (or, symmetrically, for $u_{b}=u^{\star}$ with $p=0$ independent of $u_{a}$) or for $u_{a}=u_{b}=u^{\star}$ independent of $p\in[0,1]$. In any case, one has $\tilde F(s,\tilde u^{\star})=F(s,u^{\star})$.

In case iii, the function $Q(s,\lambda,\cdot)$ is maximized for a unique value of $s_{{\rm{r}}}^{\star}=\hat s_r^{\star}(s_1)=\hat s_r^{\star}(s_2)$ independent of $\alpha$. Thus, $\tilde Q(s,\lambda,\cdot)$ is maximized when $s_{{\rm{r}}a}^{\star}$ is equal to this value with $p=1$ independent of $u_{b}$ (and, symmetrically, when $s_{{\rm{r}}b}^{\star}$ is equal to this value with $p=0$ independent of $u_{a}$) or when both $s_{{\rm{r}}a}^{\star}$ and $s_{{\rm{r}}b}^{\star}$ are equal to this value independent of $\alpha_{a}$, $\alpha_{b}$ and $p$. In any case, one has $\tilde F(s,\tilde u^{\star})=F(s,u^{\star})$, where $u^{\star}\in [0,1]\times\{s_{{\rm{r}}}^{\star}\}$.

In case iv, the function $Q(s,\lambda,\cdot)$ admits two possible maximizers. Thus, $\tilde Q(s,\lambda,\cdot)$ is maximized when $u_{a}$ is equal to one of these maximizers with $p=1$ independent of $u_{b}$, when, symmetrically, $u_{b}$ is equal to one of these maximizers with $p=0$ independent of $u_{a}$, or when $u_{a}$ and $u_{b}$ are equal to the two different maximizers independent of $p$. All these cases appear in the set-membership (\ref{chattering}).
\qed

\begin{remark}
In case 2 of Proposition \ref{prop_tilde_u}, a relaxed control $\tilde u^{\star}$ with $p\in (0,1)$
can be approximated by a high-frequency switching between non-relaxed controls  $u=(1,\hat s_{{\rm{r}}}^{\star}(s_{1}))$ and $u=(0,\hat s_{{\rm{r}}}^{\star}(s_{2}))$ (see the ``chattering control'' in \cite{BZ94}). In practice, such a high-frequency switching between the two pumps is not desired.
\end{remark}

The following Lemma will be crucial later at several places.

\begin{lemma}
\label{lemma_dot_eta}
Along any extremal trajectory, one has at almost any time
\begin{equation}
\label{dot_eta}
\dot \eta=
d\left(\frac{\gamma(s_{1})}{r}+\frac{\gamma(s_{2})}{1-r}\right)
\left(\frac{\lambda_{2}}{1-r}-\frac{\lambda_{1}}{r}\right)
+d\left(\frac{\lambda_{1}}{r^{2}}\mu(\hat s_{{\rm{r}}}^{\star}(s_{1}))+
\frac{\lambda_{2}}{(1-r)^{2}}\mu(\hat
s_{{\rm{r}}}^{\star}(s_{2}))\right)(s_{1}-s_{2}) \ .
\end{equation} 
\end{lemma}

\proof
Let us write the time derivatives of the products $\lambda_{1}\gamma(s_{1})$ and $\lambda_{2}\gamma(s_{2})$ that appear in the expression of the function $\eta$ using expressions (\ref{dynamics_convexified}), (\ref{dyn-adjoint}) and (\ref{gamma_prime}):
\[
\frac{d}{dt}\left[\lambda_{1}\gamma(s_{1})\right]=
\frac{\lambda_{1}}{r}\delta_{1}+d\gamma(s_{1})\left(\frac{\lambda_{1}}{r}-
\frac{\lambda_{2}}{1-r}\right)+d\frac{\lambda_{1}}{r}\mu(\hat s_{{\rm{r}}}^{\star}(s_{1}))(s_{2}-s_{1})
\]
where we put
\[
\delta_{1}=p\alpha_{a}\left[\mu(s_{{\rm{r}}a}^{\star})\gamma(s_{1})-\mu(\hat
  s_{{\rm{r}}}^{\star}(s_{1}))\beta(s_{1},s_{{\rm{r}}a}^{\star})\right]+
(1-p)\alpha_{b}\left[\mu(s_{{\rm{r}}b}^{\star})\gamma(s_{1})-\mu(\hat
  s_{{\rm{r}}}^{\star}(s_{1}))\beta(s_{1},s_{{\rm{r}}b}^{\star})\right].
\]
One can easily check that for any optimal control $\tilde u^{\star}$ given by Proposition \ref{prop_tilde_u}, one has $\delta_{1}=0$. Similarly, one can write
\[
\frac{d}{dt}\left[\lambda_{2}\gamma(s_{2})\right]=
\frac{\lambda_{2}}{1-r}\delta_{2}+d\gamma(s_{2})\left(\frac{\lambda_{2}}{1-r}-
\frac{\lambda_{1}}{r}\right)+d\frac{\lambda_{2}}{1-r}\mu(\hat s_{{\rm{r}}}^{\star}(s_{2}))(s_{1}-s_{2})
\]
where
\[
\delta_{2}=p(1-\alpha_{a})\left[\mu(s_{{\rm{r}}a}^{\star})\gamma(s_{2})-\mu(\hat
  s_{{\rm{r}}}^{\star}(s_{2}))\beta(s_{2},s_{{\rm{r}}a}^{\star})\right]+
(1-p)(1-\alpha_{b})\left[\mu(s_{{\rm{r}}b}^{\star})\gamma(s_{2})-\mu(\hat
  s_{{\rm{r}}}^{\star}(s_{2}))\beta(s_{2},s_{{\rm{r}}b}^{\star})\right],
\]
with $\delta_{2}=0$ for any  optimal control $\tilde u^{\star}$ given by
Proposition \ref{prop_tilde_u}.

Then, one obtains the equality (\ref{dot_eta}).
\qed

We now prove that the non-relaxed problem admits an optimal solution that is also optimal for the relaxed problem.

\begin{proposition}
\label{prop_relaxed}
The optimal trajectories for the problem with the convexified dynamics (\ref{dynamics_convexified}) are admissible optimal trajectories for the original dynamics (\ref{dynamics}). Furthermore, the optimal control $u^{\star}(\cdot)$ satisfies the following property
\[
s_{1}(t)\neq s_{2}(t) \Longrightarrow u^{\star}(t)=(1,\hat
s_{{\rm{r}}}^{\star}(s_{1})) \mbox{ or } u^{\star}(t)=(0,\hat
s_{{\rm{r}}}^{\star}(s_{2})) , \quad \mbox{for a.e. } t \in [0,t_{f}] \ .
\]
\end{proposition}

\proof
We will prove that the set of times whereby the optimal relaxed strategy generates a velocity that belongs to the convexified velocity set but not to the original velocity set has Lebesgue measure zero. For this, consider $s_{1}>s_{2}$ and $\eta=0$. Because $\gamma(\cdot)$ is increasing (see Lemma \ref{lemma3}), $\gamma(s_{1})>\gamma(s_{2})$. Additionally, $\eta=0$ implies that $\lambda_{1}/r>\lambda_{2}/(1-r)$. From equation (\ref{dot_eta}) of Lemma \ref{lemma_dot_eta}, we deduce the inequality $\dot \eta <0$ (where $\lambda_{1}$ and $\lambda_{2}$ are negative by Lemma \ref{lemma2}). Similarly, to consider $s_{2}>s_{1}$ and $\eta=0$ implies that $\dot \eta >0$. We conclude that case 2 of Proposition \ref{prop_tilde_u} can only occur at times in a set of null measure, from which the statement follows.

Now, because the optimal strategy of the convexified problem is (at almost any time) an admissible extremal for the original problem, and because the optimal time of the convexified problem is less than or equal to the optimal time of the original problem, the original problem has a solution, and it is characterized by point 1 of Proposition \ref{prop_tilde_u}. 

The last statement of the proposition follows from point 1 of Proposition \ref{prop_tilde_u}.
\qed

\section{Synthesis of the optimal strategy}
\label{sec_feedback}

According to Proposition \ref{prop_relaxed}, we can now consider optimal
trajectories of the original (non-relaxed) problem, knowing that the
optimal strategy is ``bang-bang'' except on a possible singular
arc that belongs to the diagonal set $\Delta:=\{ s \in
\Rset_{+}^{2} \mbox{ s.t. }  s_{1}=s_{2}\}$.

\begin{proposition}
\label{prop_optimal_feedback}
For $d>0$, the following feedback control drives any initial
state in $\Rset_{+}^{2}\setminus {\cal T}$ to the target ${\cal T}$ in minimal time:
\begin{equation}
\label{MRAP}
u^{\star}[s] = \left|\begin{array}{ll}
(1,\hat s_{{\rm{r}}}^{\star}(s_{1})) & \mbox{when } s_{1}>s_{2} \ , \\
(r,\hat s_{{\rm{r}}}^{\star}(s_{1}))=(r,\hat s_{{\rm{r}}}^{\star}(s_{2})) & \mbox{when } s_{1}=s_{2} \ ,\\
(0,\hat s_{{\rm{r}}}^{\star}(s_{2})) & \mbox{when } s_{1}<s_{2} \ .
\end{array}\right.
\end{equation}

\end{proposition}

\proof
From Pontryagin's Maximum Principle, a necessary optimality condition for an admissible trajectory is the existence of a solution to the adjoint system 
\begin{equation}
\label{dyn-adjoint-non-relaxed}
\left\{\begin{array}{lll}
\ds \dot\lambda_{1} & = & \ds
\lambda_{1}\frac{\alpha}{r}\mu(s_{{\rm{r}}}^{\star})
+d\left(\frac{\lambda_{1}}{r}-\frac{\lambda_{2}}{1-r}\right) \ ,\\[3mm]
\ds \dot\lambda_{2} & = & \ds 
\lambda_{2}\frac{1-\alpha}{1-r}\mu(s_{{\rm{r}}}^{\star})
+d\left(\frac{\lambda_{2}}{1-r}-\frac{\lambda_{1}}{r}\right)
\ ,
\end{array}\right.
\end{equation}
with the transversality conditions (\ref{transversality}) and where 
$u^{\star}=(\alpha,s_{{\rm{r}}}^{\star})$ maximizes the Hamiltonian
\[
H(s,\lambda,u)=-1+Q(s,\lambda,u)+d(s_{2}-s_{1})\left(\frac{\lambda_{1}}{r}-\frac{\lambda_{2}}{1-r}\right)
\]
w.r.t. $u$.\\

Consider the set
\[
I_{-}=\left\{ (s,\eta)\in(\Rset_{+}^{2}\setminus{\cal T})\times\Rset \mbox{ s.t. } s_{1}>s_{2} \mbox{ and } \eta<0\right\} \ .
\]
From expression (\ref{dot_eta}), one obtains the property
\[
s_{1}>s_{2} \mbox{ and } \eta<0 \; \Rightarrow \; \dot\eta <0
\]
using the facts that $\lambda_{i}$ ($i=1, 2$) are negative (Lemma \ref{lemma2}) and that $\gamma(\cdot)$ is increasing (Lemma \ref{lemma3}). When $\eta<0$, one has $u^{\star}=(0,\hat s_{{\rm{r}}}^{\star}(s_{2}))$ from
Lemma \ref{lemma4}, and it is possible to write
\[
\dot s_{1}-\dot s_{2} =
-\frac{d}{r(1-r)}(s_{1}-s_{2})+\frac{\gamma(s_{2})}{1-r},
\]
which shows that $s_{1}-s_{2}$ remains positive for any future
time. Thus, the set $I_{-}$ is positively invariant by the dynamics
defined by systems (\ref{reduced}) and \ref{dyn-adjoint-non-relaxed}). We deduce that the existence of a time $t<t_{f}$ such that $(s(t),\eta(t))\in I_{-}$ implies $(s(t_{f}),\eta(t_{f}))\in I_{-}$, and from the transversality condition (\ref{transversality}), one obtains $\lambda_{1}(t_{f})<\lambda_{2}(t_{f})=0$. Then, one should have $\eta(t_{f})=-\lambda_{1}(t_{f})\gamma(s_{1}(t_{f}))/r>0$, thus obtaining a contradiction. Similarly, one can show that the set
\[
I_{+}=\left\{ (s,\eta)\in(\Rset_{+}^{2}\setminus{\cal T})\times\Rset
\mbox{ s.t. } s_{1}<s_{2} \mbox{ and } \eta>0\right\}
\]
is positively invariant and that the transversality condition implies that $(s,\eta)$ never belongs to $I_{+}$ along an optimal trajectory. Because $\Delta$ is the only possible locus of a singular arc, we can form a conclusion about the optimality of (\ref{MRAP}) outside $\Delta$.\\

Now, consider the function
\[
L(s)=\frac{1}{2}(s_{1}-s_{2})^{2}
\]
and write its time derivative along an admissible trajectory $s(\cdot)$ as follows:
\[
\dot L = \langle\nabla L,\dot s\rangle=\left(-\frac{\alpha}{r}\beta(s_{1},s_{{\rm{r}}}^{\star})
+\frac{1-\alpha}{1-r}\beta(s_{2},s_{{\rm{r}}}^{\star})\right)(s_{1}-s_{2})
-\frac{2d}{r(1-r)}L \ .
\]
Along an optimal trajectory, one has
\[
\dot L+\frac{2d}{r(1-r)}L=\left|\begin{array}{ll}
\ds -\frac{\gamma(s_{1})}{r}(s_{1}-s_{2}) & \mbox{when
} s_{1}>s_{2},\\[2mm]
\ds \,\,\,\frac{\gamma(s_{2})}{1-r}(s_{1}-s_{2}) & \mbox{when
} s_{1}<s_{2},
\end{array}\right.
\]
and deduces that the inequality $\dot L+\frac{2d}{r(1-r)}L\leq 0$ is satisfied. Consequently, the set $\Delta\subset L^{-1}(0)$ is positively invariant by the optimal dynamics. On $\Delta$, the maximization of $Q(s,\lambda,\cdot)$ gives the unique $s_{{\rm{r}}}^{*}=\hat s_{{\rm{r}}}^{\star}(s_{1})=\hat s_{{\rm{r}}}^{\star}(s_{2})$ because $\lambda_{1}$, $\lambda_{2}$ are both negative (see Lemmas \ref{lemma2}, \ref{lemma3} and \ref{lemma4}). Finally, the only (non-relaxed) control that leaves $\Delta$ invariant is such that $\alpha=r$.
\qed

\begin{remark}
\label{rem_unreduced}
The feedback (\ref{MRAP}) has been proved to be optimal for the
reduced dynamics (\ref{reduced}).
In the Appendix, we prove that this feedback drives the
state of the un-reduced dynamics (\ref{slowfast}) to the target in
finite time, whatever is $\epsilon>0$. In Section \ref{sec_numerics},
we show on numerical simulations how the time to reach the
target is close from the minimal time of the reduced dynamics
when $\epsilon$ is small.
\end{remark}

\section{Study of the minimal-time function}
\label{sec_valuefunction}

Define the function
\[
T(\sigma)=\max(0,\overline T(\sigma)) \quad \mbox{with} \quad
\overline T(\sigma)=\int_{\underline s}^{\sigma}\frac{d\xi}{\gamma(\xi)}, \quad \sigma >0 \ .
\]

\begin{lemma}
\label{lemma_T}
$T(\cdot)$ is strictly concave on $[\underline s,+\infty)$.
\end{lemma}

\proof
Lemma \ref{lemma3} allows one to claim that $\overline T(\cdot)$ is twice differentiable for any $\sigma>0$ and that one has
\[
\overline T^{\prime\prime}(\sigma)=-\frac{\gamma^{\prime}(\sigma)}{\gamma(\sigma)^{2}}<0 \ ,
\quad \forall \sigma>0 \ .
\]
The function $\overline T(\cdot)$ is strictly concave on $\Rset_{+}$, and because $T(\cdot)$ coincides with $\overline T(\cdot)$ on $[\underline s,+\infty)$, we conclude that $T(\cdot)$ is strictly concave on this interval.
\qed

Let us denote the minimal-time function by $V_{d}(\cdot)$, indexed by the value of the parameter $d$:
\[
V_{d}(x)=\inf_{u(\cdot)}\left\{ t >0 \, \vert \, s(x,u,d,t)\in
{\cal T}\right\},
\]
where $s(x,u,d,\cdot)$ denotes the solution of (\ref{dynamics}) with the initial condition $s(0)=x=(x_1,x_2)$, the admissible control $u(\cdot)$ and the parameter value $d$. Lemma \ref{lemma1} ensures that these functions are well defined on $\Rset_{+}^{2}$.

\begin{proposition}\label{prop:prop4}
The value functions $V_{d}(\cdot)$ satisfy the following properties.
\begin{enumerate}[label=\roman{*}.]
\item For any $d\geq 0$, $V_{d}(\cdot)$ is Lipschitz continuous on $\Rset_{+}^{2}$.
\item For $d=0$, one has $V_{0}(x)=rT(x_{1})+(1-r)T(x_{2})$ for any $x\in \Rset_{+}^{2}$, and the feedback (\ref{MRAP}) is optimal for both relaxed and non-relaxed problems.
\end{enumerate}
\end{proposition}

\proof
On the boundary $\partial^{+}{\cal T}$ of the target that lies in the
interior of the (positively) invariant domain $\Rset_{+}^{2}$, the set
$N(\cdot)$ of unitary external normals is
\[
N(s)= \left|\begin{array}{ll}
\left\{\left(\begin{array}{c}0\\1\end{array}\right)\right\} &
\mbox{when } s_{1}<\underline s \mbox{ and } s_{2}=\underline s ,\\[4mm]
\left\{\left(\begin{array}{c}\cos \theta\\\sin
    \theta\end{array}\right)\right\}_{\theta \in [0,\pi/2]} &
\mbox{when } s_{1}=s_{2}=\underline s ,\\[4mm]
\left\{\left(\begin{array}{c}1\\0\end{array}\right)\right\} & \mbox{when } s_{1}=\underline s \mbox{ and } s_{2}<\underline s \ .
\end{array}\right.
\]
At any $s \in \partial^{+}{\cal T}$, one has
\[
\inf_{u \in U(s)}\inf_{\nu \in N(s)} \langle F(s,u)+dG(s),\nu\rangle \leq \inf_{u \in U(s)}\inf_{\nu \in N(s)} \langle
F(s,u),\nu\rangle=-\gamma(\underline s)<0 \ .
\]
Furthermore, the maps
\[
s \mapsto F(s,u)+dG(s)
\]
are Lipschitz continuous w.r.t.~$s \in \Rset_{+}^{2}$ uniformly in $u$. According to \cite[Sect 1. and 4, Chap. IV]{BC97}, the target satisfies then the {\em small time locally controllable} property, and the value functions $V_{d}(\cdot)$ are Lipschitz continuous on $\Rset_{+}^{2}$.\\

When $d=0$, the feedback (\ref{MRAP}) provides the following dynamics
\[
\dot s = \left|\begin{array}{rl}
-\frac{1}{r}\left[\begin{array}{c}\gamma(s_{1})\\0\end{array}\right]
& \mbox{when } s_{1}>\max(s_{2},\underline s) \ ,\\[4mm]
-\left[\begin{array}{c}\gamma(s_{1})\\\gamma(s_{2})\end{array}\right]
& \mbox{when } s_{1}=s_{2}>\underline s \ , \\[4mm]
-\frac{1}{1-r}\left[\begin{array}{c}0\\\gamma(s_{2})\end{array}\right]
& \mbox{when } s_{2}>\max(s_{1},\underline s) \ ,
\end{array}\right.
\]
and one can explicitly calculate the time to go to the target for any initial condition $x \in \Rset_{+}^{2}$, which we denote as $W_{0}(x)$:
\[
W_{0}(x)=\left|\begin{array}{ll}
\ds -r\int_{x_{1}}^{\max(x_{2},\underline s)}\frac{ds_{1}}{\gamma(s_{1})}
-\int_{\max(x_{2},\underline s)}^{\underline s}\frac{ds_{1}}{\gamma(s_{1})}
& \mbox{when } x_{1}\geq x_{2} \ , \\[4mm]
\ds -(1-r)\int_{x_{2}}^{\max(x_{1},\underline s)}\frac{ds_{2}}{\gamma(s_{2})}
-\int_{\max(x_{1},\underline s)}^{\underline s}\frac{ds_{2}}{\gamma(s_{2})}
& \mbox{when } x_{1}\leq x_{2} \ .
\end{array}\right.
\]
One can check that $W_{0}$ is Lipschitz continuous and that it can be written as $W_{0}(x)=rT(x_{1})+(1-r)T(x_{2})$. We now show that $W_{0}$ is a viscosity solution of the Hamilton-Jacobi-Bellman equation associated to the relaxed problem
\begin{equation}
\label{HJB0}
{\cal H}(x,\nabla W_{0}(x))=-1+\max_{(u_{a},u_{b},p)\in\tilde U(x)} pQ(x,-\nabla
W_{0}(x),u_{a})+(1-p)Q(x,-\nabla W_{0}(x),u_{b})=0, \quad x \notin {\cal T},
\end{equation}
(where $Q$ is defined in (\ref{defQ})) with the boundary condition
\begin{equation}
\label{BC}
W_{0}(x)=0, \quad x \in {\cal T} \ .
\end{equation}
Consider the $C^{1}$ functions
\[
\overline{W}_{0,1}(x) = r\overline T(x_{1}) \quad , \quad
\overline{W}_{0,2}(x) = (1-r)\overline T(x_{2}) \quad \mbox{and} \quad 
\overline{W}_{0}(x) = \overline{W}_{0,1}(x)+\overline{W}_{0,2}(x)
\]
defined on $\Rset_{+}^{2}$.
One has
\[
\nabla \overline W_{0,1}(x)=\left[\begin{array}{c}
\ds \frac{r}{\gamma(x_{1})}\\[4mm]
0
\end{array}\right] \quad \mbox{and} \quad
\nabla \overline W_{0,2}(x)=\left[\begin{array}{c}
0\\[4mm]
\ds \frac{1-r}{\gamma(x_{2})}
\end{array}\right] 
\]
which are non-negative vectors. One can then use Lemma \ref{lemma4} to obtain the property
\[
{\cal H}(x,\nabla \overline W_{0,1}(x))=
{\cal H}(x,\nabla \overline W_{0,2}(x))=
{\cal H}(x,\nabla \overline W_{0}(x))=0, \quad x \in \Rset_{+}^{2} \ ,
\]
which shows that $\overline W_{0,1}$, $\overline W_{0,2}$ and $\overline W_{0}$ are solutions of (\ref{HJB0}) in the classical sense.

At $x \notin {\cal T}$ with $x_{i}\neq \underline s$ ($i=1,2$), $W_{0}$ is $C^{1}$ and locally coincides with $\overline W_{0}$. Then, it satisfies equation (\ref{HJB0}) in the classical sense.

At $x \notin {\cal T}$ with $x_{1}= \underline s$ or $x_{2}=\underline s$, $W_{0}$ is not differentiable but locally coincides with $\max(\overline W_{0},\overline W_{0,2})$ or $\max(\overline W_{0},\overline W_{0,1})$. From the properties of viscosity solutions (see, for instance, \cite[Prop 2.1, Chap. II]{BC97}), one must simply check that $W_{0}$ is a super-solution of (\ref{HJB0}). At such points, the Fr\'echet
sub-differential of $W_{0}$ is
\[
\partial^{-}W_{0}(x) = \left|\begin{array}{ll}
\ds \left[0,\frac{r}{\gamma(\underline
  s)}\right]\times\left\{\frac{1-r}{\gamma(x_{2})}\right\}
& \mbox{when } x_{1}=\underline s \ ,\\[4mm]
\ds \left\{\frac{r}{\gamma(x_{1})}\right\}\times\left[0,\frac{1-r}{\gamma(\underline
  s)}\right]
& \mbox{when } x_{2}=\underline s \ .
\end{array}\right.
\]
Because any sub-gradient $\delta^{-}\in \partial^{-}W_{0}(x)$ is a non-negative vector, one can again use Lemma \ref{lemma4} and obtain
\[
{\cal H}(x,\delta^{-})=0, \quad \forall \delta^{-}
\in \partial^{-}W_{0}(x) \ ,
\]
which proves that $W_{0}$ is a viscosity solution of \eqref{HJB0}. Moreover, $W_{0}$ satisfies the boundary condition (\ref{BC}). Finally, we use the characterization of the minimal-time function as the unique viscosity solution of (\ref{HJB0}) in the class of Lipschitz continuous functions with boundary conditions (\ref{BC}) (see \cite[Th. 2.6, Chap IV]{BC97}) to conclude that $W_{0}$ is the value function of the relaxed problem. Because the time $W_{0}(x)$ to reach the target from an initial condition $x \notin {\cal T}$ is obtained with the non-relaxed control (\ref{MRAP}), we also deduce that $V_{0}$ and $W_{0}$ are equal.
\qed

\begin{remark}
In the case $d=0$, the control given by \eqref{MRAP} is optimal but not the unique solution of the problem. Indeed, in Proposition \ref{prop:prop4}, we proved that $V_0(\cdot)$ is the unique viscosity solution to equation \eqref{HJB0}, where one of the possible maximizers of the Hamiltonian given in \eqref{HJB0} is given by \eqref{MRAP}, but on the set $(\underline s,\infty)^2\setminus \Delta$ there are more choices for $u$; for instance,
\begin{equation*}
u^{\star}[s] = \left|\begin{array}{ll}
(1,\hat s_{{\rm{r}}}^{\star}(s_{1})) & \mbox{when } s_{2}\leq\underline s<s_{1} \ , \\
(0,\hat s_{{\rm{r}}}^{\star}(s_{2})) & \mbox{when } s_{1}>s_{2}>\underline s \ , \\
(r,\hat s_{{\rm{r}}}^{\star}(s_{1}))=(r,\hat s_{{\rm{r}}}^{\star}(s_{2})) & \mbox{when } s_{1}=s_{2} \ ,\\
(1,\hat s_{{\rm{r}}}^{\star}(s_{1})) & \mbox{when } \underline s<s_{1}<s_{2} \ ,\\
(0,\hat s_{{\rm{r}}}^{\star}(s_{2})) & \mbox{when } s_{1}\leq \underline s < s_{2} \ 
\end{array}\right.
\end{equation*}
satisfies \eqref{HJB0}.
\end{remark}

\begin{proposition}
The functions $V_{d}(\cdot)$ satisfy the following properties:
\begin{enumerate}[label=\roman{*}.]
\item $V_{d}(x)=T(x_{1})=T(x_{2})$ for any $x \in \Delta$ and $d\geq 0$,
\item $V_{\infty}(x)=\lim_{d\to+\infty}V_{d}(x)=T(rx_{1}+(1-r)x_{2})$ for any $x\in\Rset_{+}^{2}$, and
\item $d \mapsto V_{d}(x)$ is increasing for any $x \in (\underline s,+\infty)^{2}\setminus \Delta$.
\end{enumerate}
\end{proposition}

\proof
Consider an initial condition $x$ in $\Delta\setminus{\cal T}$. The optimal synthesis given in Proposition \ref{prop_optimal_feedback} shows that the set $\Delta$ is invariant by the optimal flow and that the dynamics on $\Delta$ are
\[
\dot s_{i}=-\gamma(s_{i}), \quad i=1, 2
\]
independent of $d$. We then conclude that $V_{d}(x)=T(x_{i})$ for $i=1,
2$.\\

Consider $d>0$ and $x\notin\Delta\cup{\cal T}$. Denote for simplicity $s(\cdot)$ as the solution $s(x,u^{\star},d,\cdot)$ with the feedback control $u^{\star}$ given in Proposition \ref{prop_optimal_feedback}, and $t_{f}=V_{d}(x)$. Define $t_{\Delta}$ as the first time $t$ such that $s(t)\in \Delta$ (here, we allow the solution $s(\cdot)$ to possibly enter the target ${\cal T}$ before reaching $\Delta$). 

From equation (\ref{dynamics}) with control (\ref{MRAP}), one can easily check that the following inequalities are satisfied
\[
\begin{array}{l}
x_{1}>x_{2} \Rightarrow x_{1} > s_{1}(t) \geq s_{2}(t) > x_{2}, \quad
\forall t\in [0,t_{\Delta}] \ ,\\[2mm]
x_{1}<x_{2} \Rightarrow x_{1} < s_{1}(t) \leq s_{2}(t) < x_{2}, \quad
\forall t\in [0,t_{\Delta}] \ .
\end{array}
\]
Then, because the function $\gamma(\cdot)$ is increasing (Lemma \ref{lemma3}),
one can write, if the state $s$ has not yet reached $\Delta$,
\begin{equation}
\label{bounds_d}
-\frac{d}{r(1-r)}|s_{1}-s_{2}|-M_{+}\leq \frac{d}{dt}|s_{1}-s_{2}|
\leq -\frac{d}{r(1-r)}|s_{1}-s_{2}|-M_{-}
\end{equation}
with $M_{-}=\min(\gamma(x_{2})/r,\gamma(x_{1})/(1-r))$ and $M_{+}=\max(\gamma(x_{1})/r,\gamma(x_{2})/(1-r))$. Then, we obtain an upper bound on the time $t_{\Delta}$
\begin{equation}
\label{bounds_t}
t_{\Delta} \leq
\frac{r(1-r)}{d}\log\left(1+d\frac{|x_{1}-x_{2}|}{M_{-}r(1-r)}\right)
\end{equation}
which tends to zero when $d$ tends to infinity. From (\ref{bounds_d}), we can also write
\[
|x_{1}-x_{2}|-M_{+}t_{\Delta}
\leq \frac{d}{r(1-r)}\int_{0}^{t_{\Delta}}|s_{2}(\tau)-s_{1}(\tau)|d\tau
\leq |x_{1}-x_{2}|-M_{-}t_{\Delta}
\]
and finally, one obtains from (\ref{dynamics}) the following bounds on $s_{i}(t_{\Delta})$ ($i=1,2$):
\begin{equation}
\label{bounds_s}
rx_{1}+(1-r)x_{2}-\max(r,(1-r))M_{+}t_{\Delta} \leq
s_{i}(t_{\Delta})\leq
rx_{1}+(1-r)x_{2}-\min(r,(1-r))M_{-}t_{\Delta} \ .
\end{equation}
Therefore, $s_{1}(t_{\Delta})=s_{2}(t_{\Delta})$ converges to $rx_{1}+(1-r)x_{2}$ 
when $d$ tends to $+\infty$. Furthermore, one has
\[
\begin{array}{ll}
t_{f} = t_{\Delta} + T(s(t_{\Delta})) & \mbox{when } s(t_{\Delta})
  \notin {\cal T} \ ,\\
t_{f} < t_{\Delta} & \mbox{when } s(t_{\Delta}) \in {\cal T} \ .
\end{array}
\]
Because $t_{\Delta}\to 0$ and because $T(\cdot)$ is continuous with $T(rx_{1}+(1-r)x_{2})=0$ when $rx_{1}+(1-r)x_{2}\leq \underline s$, we obtain the convergence
\[
V_{\infty}(x)=\lim_{d\to+\infty} V_{d}(x)=T(rx_{1}+(1-r)x_{2}) \ .
\]

\bigskip

Now, consider the domain ${\cal D}_{+}=\{ s \in \Rset_{+}^{2} \, \vert \, s_{1}\geq s_{2}>\underline s\}$, and let us show that any trajectory of the optimal flow leaves ${\cal D}_{+}$ at $(\underline s,\underline s)$ with the help of this simple argumentation on the boundaries of the domain:
\[
\begin{array}{l}
\ds s_{2}=\underline s \; \Rightarrow \;
\dot s_{2}=\frac{d}{1-r}(s_{1}-\underline s)\geq 0 \ ,\\[3mm]
\ds s_{1}=s_{2} \; \Rightarrow \;
\dot s_{1}=\dot s_{2} \ .
\end{array}
\]
It is convenient to consider the variable $\tilde s=rs_{1}+(1-r)s_{2}$, whose optimal dynamics in ${\cal D}_{+}$ are simply
\begin{equation}
\label{dyn_s_tilde}
\dot{\tilde s}(t) = -\gamma(s_{1}(t)) \ , \quad t \in [0,t_{f}] \ .
\end{equation}
Because $\tilde s(\cdot)$ is strictly decreasing with time, an optimal trajectory in ${\cal D}_{+}$ can be parameterized by the fictitious time
\begin{equation}
\label{fictious_time}
\tau(t)=rx_{1}+(1-r)x_{2}-\tilde s(t) \ , \quad t \in [0,t_{f}]
\end{equation}
(where $x$ is an initial condition in ${\cal D}_{+}$). The variable $s_{1}(\cdot)$ is then a solution of the scalar non-autonomous dynamics
\[
\frac{ds_{1}}{d\tau}=f_{d}(\tau,s_{1})=
\left|\begin{array}{ll}
\ds -\frac{1}{r}-d\frac{s_{1}+\tau-(rx_{1}+(1-r)x_{2})}{r(1-r)\gamma(s_{1})}
& \mbox{when } s_{1}+\tau>rx_{1}+(1-r)x_{2} \ , \\[3mm]
-1 & \mbox{when } s_{1}+\tau=rx_{1}+(1-r)x_{2} \ ,
\end{array}\right.
\]
with the terminal fictitious time
\[
\tau_{f}=rx_{1}+(1-r)x_{2}-\underline s \ .
\]
Notice that $\tau_{f}$ is independent of $d$. One then deduces the inequalities
\[
d_{1} > d_{2} \mbox{ and } s_{1}+\tau>rx_{1}+(1-r)x_{2}
\Longrightarrow f_{d_{1}}(\tau,s_{1})<f_{d_{2}}(\tau,s_{1}) 
\]
and thus,
\begin{equation}
\label{ineq_s1}
d_{1} > d_{2} \mbox{ and } x \in {\cal D}_{+}\setminus\Delta
\Longrightarrow
s_{1}(x,u^{\star},d_{1},\tau)<s_{1}(x,u^{\star},d_{2},\tau) \ , 
\; \forall \tau \in [0,\tau_{f}] \ .
\end{equation}
Finally, from equations (\ref{dyn_s_tilde}) and (\ref{fictious_time}), the time to reach the target can be expressed as
\begin{equation}
\label{tf}
t_{f}=\int_{0}^{\tau_{f}} \frac{d\tau}{\gamma(s_{1}(\tau))} \ .
\end{equation}
Because the function $\gamma(\cdot)$ is increasing and because $\tau_{f}$ is independent of $d$, one can conclude from (\ref{ineq_s1}) and (\ref{tf}) that
\[
d_{1} > d_{2} \mbox{ and } x \in {\cal D}_{+}\setminus\Delta \Longrightarrow V_{d_{1}}(x)>V_{d_{2}}(x) \ .
\]
The case of initial conditions in ${\cal D}_{-}\setminus\Delta$, with ${\cal D}_{-}=\{ s \in \Rset_{+}^{2} \, \vert \, s_{2}\geq s_{1}>\underline s\}$, is symmetric.
\qed

\begin{remark}\label{remark4}
The tightness $V_{\infty}-V_{0}$ of the bounds on the value function
$V_{d}$ on $(\underline s,+\infty)^{2}\setminus \Delta$ is related to
the concavity of the function $T(\cdot)$ on $(\underline s,+\infty)$
(the less the concavity $\max_{\sigma\in[\underline
    s,+\infty)}|\overline T^{\prime\prime}(\sigma)|$ is, the tighter the bounds are).

The bounds $V_{0}\leq V_{d}<V_{\infty}$ that are satisfied on the set $(\underline s,+\infty)^{2}\setminus \Delta$ are not necessarily satisfied outside this set: for $x$ outside the target but such that $rx_{1}+(1-r)x_{2}<\underline s$, one has $V_{\infty}(x)=0$ and $V_{0}(x)>0$. Therefore, we conclude that a large diffusion negatively impacts the time to treat the resource when both zones are initially polluted; however, when one of the two zones is initially under the pollution threshold, a large diffusion could positively impact the duration of the treatment.
\end{remark}

\section{Numerical illustrations}
\label{sec_numerics}
\begin{figure}[h!]
\begin{center}
\includegraphics[scale=0.35]{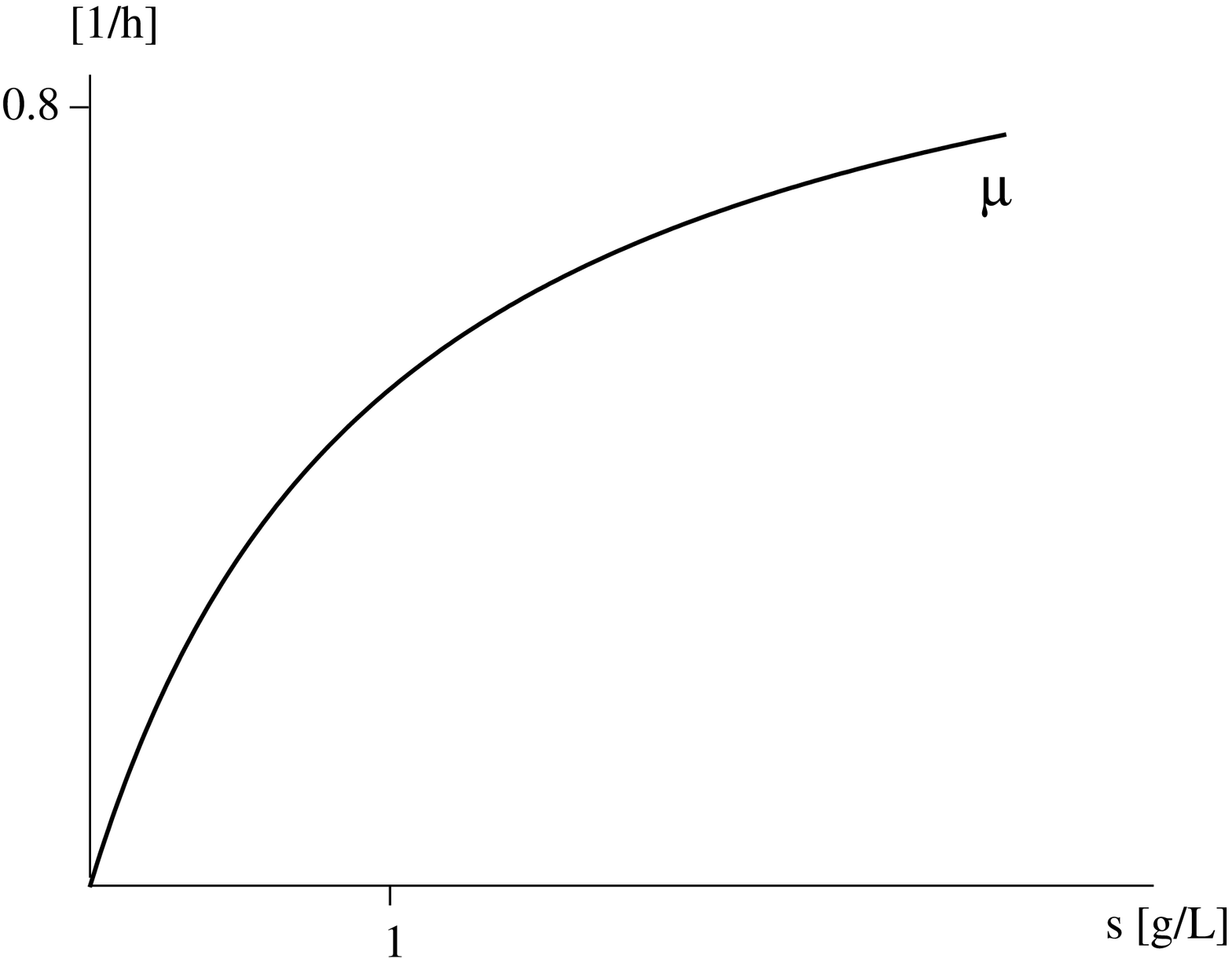}
\hspace{5mm}
\includegraphics[scale=0.35]{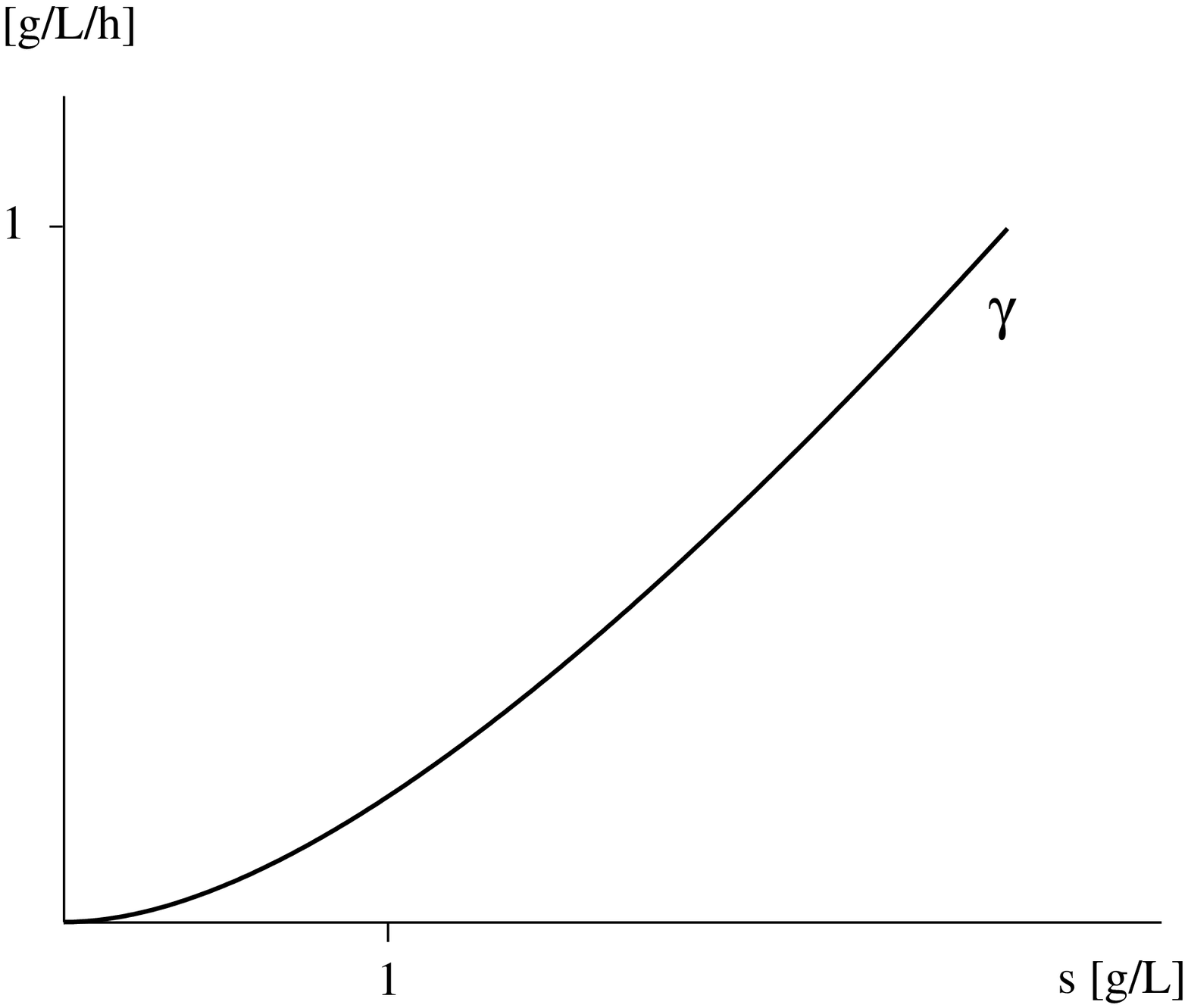}
\caption{Graphs of $\mu(\cdot)$ and corresponding $\gamma(\cdot)$.\label{fig2}}
\end{center}
\end{figure}

\begin{figure}[h!]
\begin{center}
\includegraphics[scale=0.4]{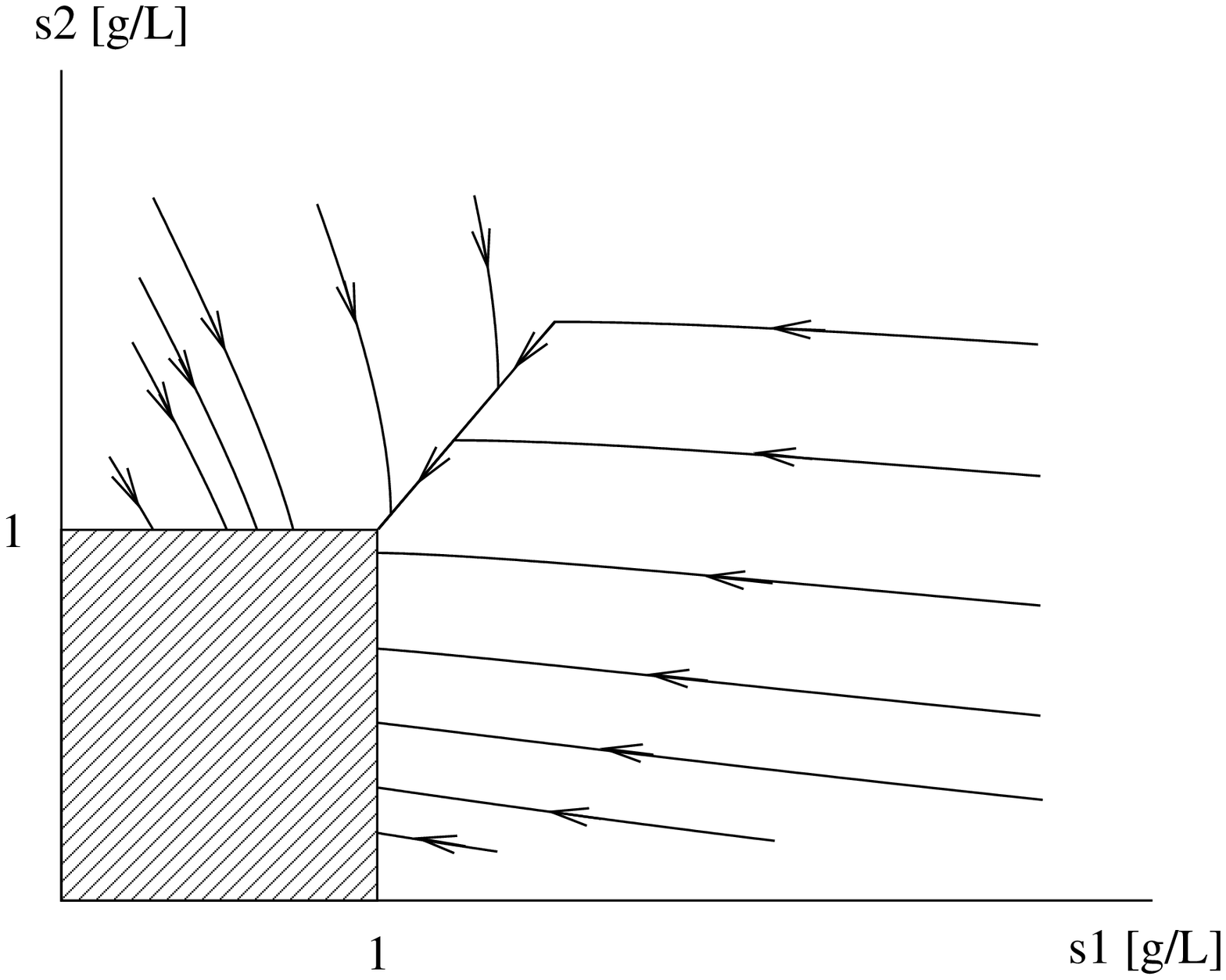}
\hspace{5mm}
\includegraphics[scale=0.4]{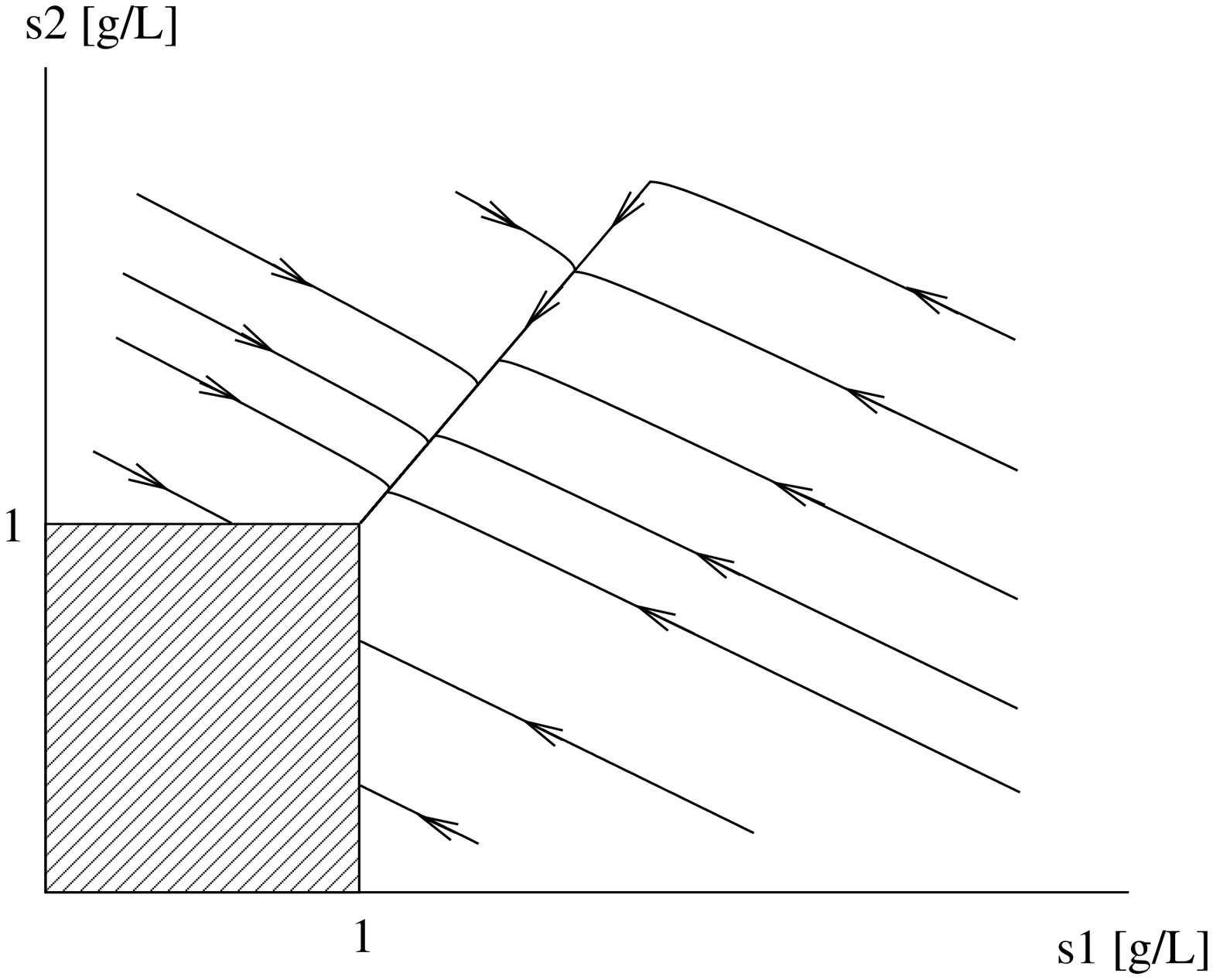}
\caption{Optimal paths for $d=0.1[h^{-1}]$ (left) and $d=10[h^{-1}]$ (right)
  with $r=0.3$.\label{fig3}}
\end{center}
\end{figure}

\begin{figure}[h!]
\begin{center}
\includegraphics[scale=0.4]{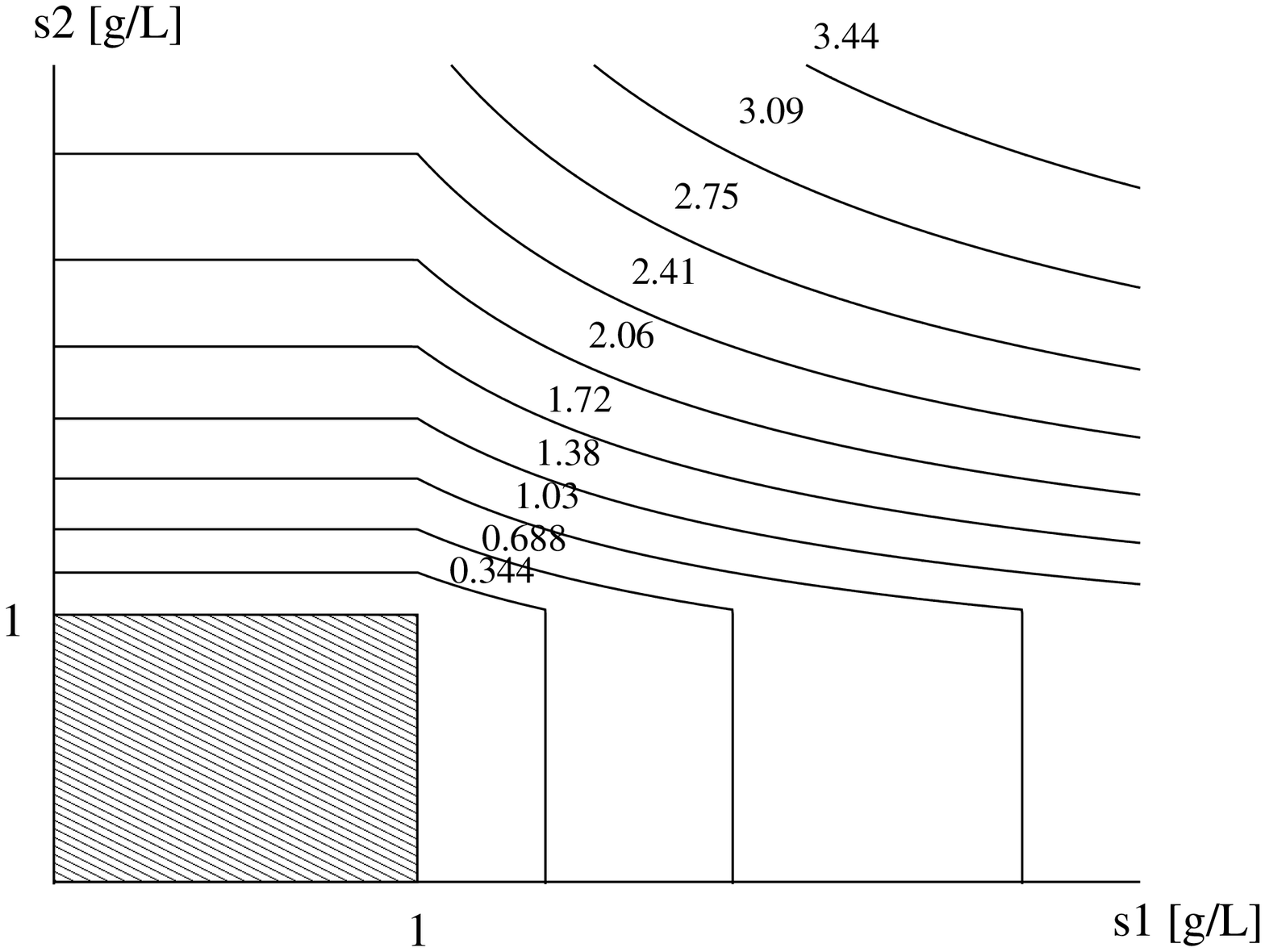}
\hspace{5mm}
\includegraphics[scale=0.4]{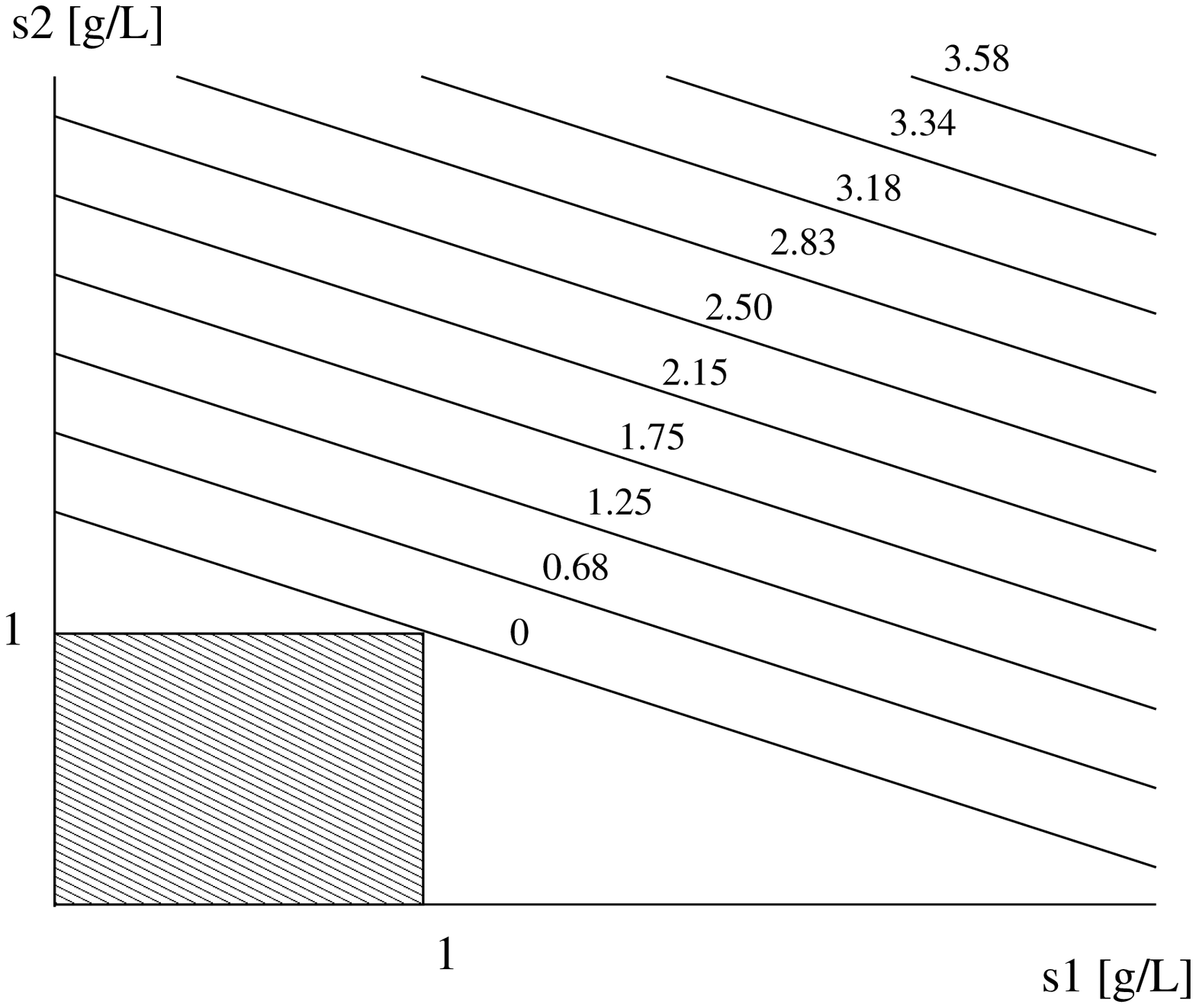}
\caption{Level sets (in hours) of $V_{0}$ (left) and $V_{\infty}$ (right)
  for $r=0.3$.\label{fig4}}
\end{center}
\end{figure}

We consider the Monod (or Michaelis-Menten) growth function, which is 
quite popular in bio-processes and which satisfies Assumption \ref{H0}:
\[
\mu(s)=\mu_{\max}\frac{s}{K_{s}+s},
\]
with the parameters $\mu_{\max}=1 [h^{-1}]$ and
$K_{s}=1 [gL^{-1}]$. The corresponding function
$\gamma(\cdot)$ is depicted in Fig. \ref{fig2}. The threshold that
defines the target has been chosen as $\underline s=1[gL^{-1}]$.
Several optimal trajectories in the phase portrait are drawn in Fig. \ref{fig3} for small and large values of the parameter $d$.
Finally, level sets of the value functions $V_{0}$ and $V_{\infty}$ are represented in Fig. \ref{fig4}.
One can make the following observations concerning the influence of
the diffusion on the treatment duration, that we consider to be
valuable from a practical viewpoint.
\begin{itemize}
\item When pollution is homogeneous, the best is to maintain it
  homogeneous, and the treatment time is then independent of the diffusion.
\item A high diffusion is favorable for having fast treatments when 
initial concentrations are strongly different for the two zones.
Typically, when the pollutant concentration is below the threshold in
one patch, a high diffusion can reduce significantly the treatment time
compared to a small diffusion.
\item When initial concentrations in the two patches are close, a
  small diffusion leads to faster treatment than a large diffusion.
\end{itemize}

\begin{figure}[h!]
\begin{center}
\includegraphics[scale=0.5]{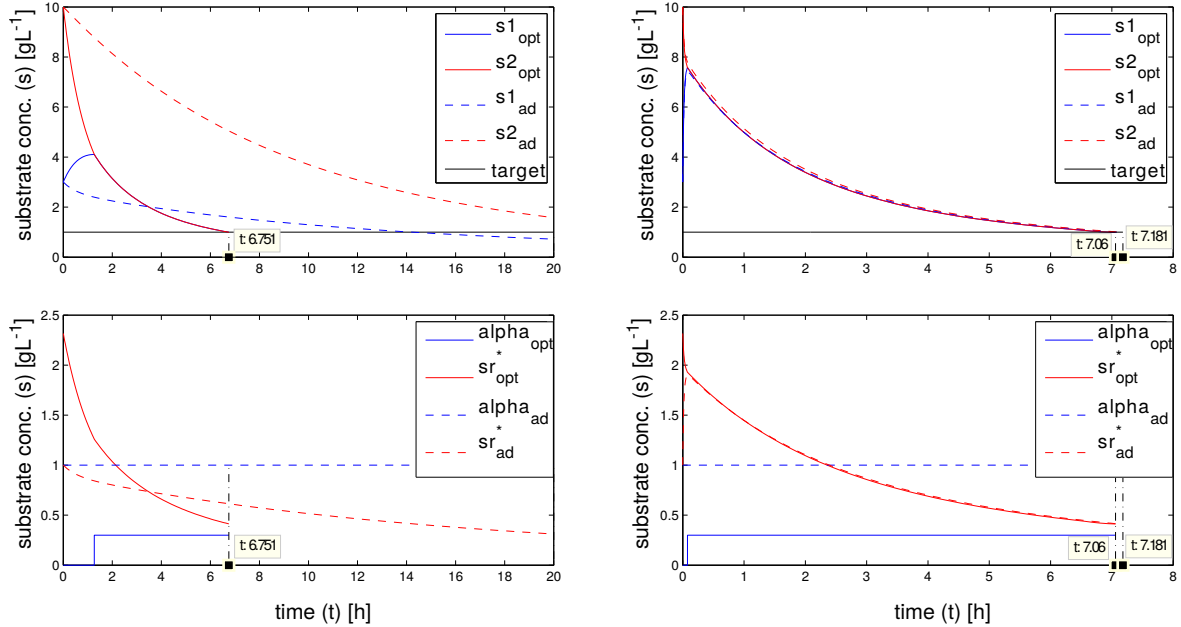}
\caption{Trajectories and controls generated by the
  {\em two-pumps} and {\em one-pump} optimal feedback, for $r=0.3$,
  $\underline s=1 [gL^{-1}]$ and $s(0)=(3,10) [gL^{-1}]$.
 On the left $d=0.1 [h^{-1}]$, and on the right $d=10[h^{-1}]$.\label{fig5}}
\end{center}
\end{figure}

For various initial condition $s(0)$, we have also performed numerical
comparisons of the minimal time $V_{d}(s(0))$ given by the feedback
strategy \eqref{MRAP} against two other non-optimal control strategies:
\begin{enumerate}
\item the {\em best constant control} that gives the smallest time
  $T^\star_{cst}$ to reach the target among constant controls,
\item the optimal {\em one-pump} feedback strategy
obtained in the former work \cite{GRRR12}. This last control strategy
considers that only one patch can be treated (that we called the
``active zone''). The problem amounts then
to consider the same dynamics (\ref{reduced}) but one seeks the
feedback $s_{{\rm{r}}}^\star(\cdot)$ that gives the minimal time
$T^\star_{one}$ when $\alpha$ is imposed to be
constantly equal to $1$ (or $0$ depending which patch is treated). 
In \cite{GRRR12}, it has been proved that the feedback $s_{1} \mapsto
\hat s_{{\rm{r}}}^\star(s_{1})$ is optimal.
\end{enumerate}
\begin{table}[ht!]
  \begin{center}
  \renewcommand{\arraystretch}{1.2}
  \begin{tabular}{||c||c|c|c|c|c|c||}
    \hline\hline
        & \multicolumn{2}{c|}{$\mathbf{V}_{d}$} &
\multicolumn{2}{c|}{$\mathbf{T_{cst}^{\star}}$} &
\multicolumn{2}{c||}{$\mathbf{T_{one}^{\star}}$} \\
    \cline{2-7}
    & $d=0.1$ & $d=10$ & $d=0.1$ & $d=10$ & $d=0.1$ & $d=10$\\
    \hline\hline
    $s(0)=(1.5,0)$ & 0.42   & 0.01   & 0.42   & 0.01   & 0.42    & 0.01
\\
    \hfill  Increase: &  &  & (+ 1.45 \%)& (+ 0.00 \%)  & (+ 0.00 \%) & (+
0.00 \%)\\ \hline
    $s(0)=(3,0)$   & 1.01   & 0.06   & 1.05   & 0.06   & 1.01   & 0.06   \\
    \hfill  Increase: &  &  & (+ 3.90 \%)& (+ 0.85 \%)  & (+ 0.00 \%)
& (+ 0.00 \%)\\ \hline
    $s(0)=(4,0.5)$  & 1.33   & 2.17   & 1.39   & 2.23   & 1.37   & 2.21
\\
    \hfill  Increase: &  &  & (+ 4.68 \%)& (+ 2.62 \%)  & (+ 2.73 \%) & (+
1.55 \%)\\ \hline
    $s(0)=(4,1.5)$  & 3.20   & 3.65   & 3.67   & 3.75   & 8.27   & 3.72
\\
    \hfill  Increase: &  &  & (+ 14.76 \%)& (+ 2.58 \%)  & (+ 158.27 \%) &
(+ 1.91 \%)\\ \hline
    $s(0)=(4,4)$    & 5.45   &  5.45  & 5.74   & 5.71   & 18.25   & 5.53
\\
    \hfill  Increase: &  &  & (+ 5.43 \%)& (+ 4.90 \%)  & (+ 235.01 \%) &
(+ 1.59 \%)\\ \hline \hline
  \end{tabular}
\vspace{3mm}

  \begin{tabular}{||c||c|c|c|c|c|c||}
    \hline\hline
       & \multicolumn{2}{c|}{$\mathbf{V}_{d}$} &
\multicolumn{2}{c|}{$\mathbf{T_{cst}^{\star}}$} &
\multicolumn{2}{c||}{$\mathbf{T_{one}^{\star}}$} \\
    \cline{2-7}
    & $d=0.1$ & $d=10$ & $d=0.1$ & $d=10$ & $d=0.1$ & $d=10$\\
    \hline\hline
    $s(0)=(1.5,0)$ & 25.95   & 34.12   & 38.65   & 38.81   & 34.03    &
34.14   \\
    \hfill  Increase: &  &  & (+ 48.93 \%)& (+ 13.74 \%)  & (+ 31.14 \%) &
(+ 0.05 \%)\\ \hline
    $s(0)=(3,0)$   & 32.91   & 39.91   & 50.08   & 50.12   & 45.89   &
40.15   \\
    \hfill  Increase: &  &  & (+ 52.18 \%)& (+ 25.58 \%)  & (+ 39.45 \%) &
(+ 0.60 \%)\\ \hline
    $s(0)=(4,0.5)$  & 41.08   & 42.86   & 58.65   & 58.02   & 61.51   &
42.94   \\
    \hfill  Increase: &  &  & (+ 42.77 \%)& (+ 35.37 \%)  & (+ 49.74 \%) &
(+ 0.1 \%)\\ \hline
    $s(0)=(4,1.5)$  & 43.69   & 44.37   & 63.59   & 63.28  & 70.81   &
44.49   \\
    \hfill  Increase: &  &  & (+ 45.57 \%)& (+ 42.61 \%)  & (+ 62.08 \%) &
(+ 0.27 \%)\\ \hline
    $s(0)=(4,4)$    & 45.94 &  45.94  & 71.67 & 71.04  & 81.58 &  46.17  \\
    \hfill  Increase: &  &  & (+ 56.02 \%)& (+ 54.64 \%)  & (+ 77.60 \%) &
(+ 0.51 \%)\\ \hline \hline
  \end{tabular}
  \caption{\label{tables} Time comparisons (in hours) for
      $r=0.3$ and target value $\underline s = 1$ (top), $0.1$
      (bottom) $[gL^{-1}]$ (initial condition $s(0)$ and
    diffusion parameter $d$ are given in $gL^{-1}$ and $h^{-1}$,
    respectively).}
\end{center}
\end{table}
The results presented in Table \ref{tables}
show first that the benefit of using the optimal feedback strategy
over the other strategies increases with the level of initial
pollution.
The simulations also
demonstrate the gain of using two pumps instead of one: for large
concentrations of pollutant at initial time, one can see on the tables that a
constant two-pumps strategy can be even better that the optimal
feedback strategy restricted to the use of one pump only.
This kind of situations typically occurs when diffusion is low and the time
required by the optimal strategy for using simultaneously the two pumps
is large compared to the overall duration. This is particularly
noticeable when the initial pollution is homogeneous and the use of two pumps
allows to maintain the levels of concentrations equal in both patches.
We conclude that, for small diffusion, treating only one patch without
the possibility to allocate the treatment in both patches could be
quite penalizing. Figure \ref{fig5}
illustrates the time history of the two feedback controllers.

Furthermore, the Table \ref{tables} illustrates the
effect of diffusion on the treatment times. 
One can first notice that the relative effect of the
diffusion parameter $d$ on the optimal time $V_{d}$ is decreasing with
the threshold $\underline s$. This can be explained by the fact that
the proportion of the time spent on the set $s_{1}=s_{2}$, that is
independent of the parameter $d$, is larger when one begins further
away from the target.
One can also see that the differences between strategies
decrease when the diffusion increases. Intuitively, a high diffusion
makes the resource behave quickly close to a perfectly mixed resource with
one patch, leading consequently to less benefit of using more than one
pump. 
Nevertheless, one can see that considering feedback controls
remain quite efficient compared to constant ones when initial
pollution is high.\\

\begin{figure}[h!]
\begin{center}
\includegraphics[scale=0.4]{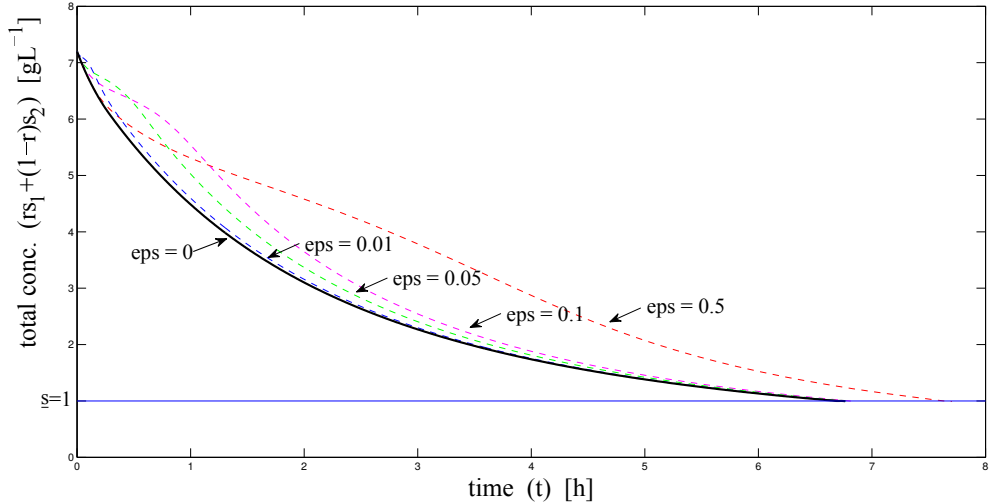}
\caption{Total pollutant concentration in the resource of
  the full dynamics \eqref{slowfast} with the strategy \eqref{MRAP}, for different values of $\epsilon$.\label{fig6}}
\end{center}
\end{figure}

Finally, we illustrate on Fig. \ref{fig6} the effect of approximating
the original dynamics \eqref{slowfast} by the reduced one
\eqref{reduced}, when applying the feedback \eqref{MRAP}.
As proven in the Appendix, the feedback \eqref{MRAP} drives the state
to the target in finite time for any $\epsilon>0$.

\section{Conclusion}
In this work, we have shown that although the velocity set of the control problem is not convex, there exists an optimal solution with ordinary controls that is also optimal among relaxed controls. The optimal strategy consists in the most rapid approach to the homogenized concentration of
pollutant in both patches. For the particular case of null diffusion, the most rapid approach path is not the unique solution of the problem. This optimal state-feedback has some interesting features for the practitioners and controllers: 
\begin{enumerate}
\item it does not require knowledge of the diffusion parameter $D$ to be implemented, and
\item if the ratio $r$ of the volumes of the two patches is not known,
  the optimal trajectory can be approximated by a regularization of
  the bang-bang control about the neighborhood of the set
  $s_{1}=s_{2}$ that keeps the trajectory in this neighborhood.
\end{enumerate}
Furthermore, is has been shown in simulations that the
    benefit of using two pumps instead of one can be significant when
    the diffusion is low.
We have also proposed explicit bounds on the minimal-time function, characterizing the extreme cases $d=0$ and $d=+\infty$. We have shown that a large diffusion rate increases the treatment time when the pollution concentration is above the desired threshold in both zones, while in contrast, it can be beneficial when the concentration in one of the two zones is below the desired threshold. This remarkable feature could serve practitioners in the choice of pump positioning in an originally clean water resource that is suddenly affected by a local pollution. Such an investigation could be the matter of future work.

\section*{Acknowledgments}

This work was developed in the context of the DYMECOS 2 INRIA Associated team and of project BIONATURE of CIRIC INRIA CHILE, and it was partially supported by CONICYT grant REDES 130067. The first and third authors were also supported by CONICYT-Chile under ACT project 10336, FONDECYT 1110888, BASAL project (Centro de Modelamiento Matem\'atico, Universidad de Chile), MathAmsud N\textdegree 15MATH-02, and Beca Doctorado Nacional Convocatoria 2013 folio 21130840 CONICYT-CHILE. The third author acknowledges the support of Departamento de Postgrado y Post\'itulo de la Vicerrector\'ia de Asuntos Acad\'emicos (Universidad de Chile) and Institut Fran\c{c}ais (Ambassade de France au Chili).

The authors are also grateful to T. Bayen, J. F. Bonnans, P. Gajardo, J. Harmand, and A. Rousseau for fruitful discussions and insightful ideas.

\section*{Appendix}

\begin{proposition}
\label{approximation}
For any $\epsilon>0$, the feedback strategy (\ref{MRAP}) applied to
the full dynamics (\ref{slowfast}) with $x_{{\rm{r}}}(0)>0$ drives the state
to the target in finite time.
\end{proposition}

\proof

Without any loss of generality, we assume that $s_{1}(0)\geq s_{2}(0)$
(the proof is similar when $s_{1}(0)\leq s_{2}(0)$).

If $s_{1}(0)> s_{2}(0)$, we prove that $s_{1}=s_{2}$ is reached in
finite time. If not one, one should have $s_{1}(t)> s_{2}(t)$ with
$s_{1}(t)\geq \underline s$ for any $t>0$. This implies to have
$\alpha^\star(t)=1$ and $s_{{\rm{r}}}^\star(t)=\hat s_{{\rm{r}}}^\star(s_{1}(t))$ at any
time $t>0$ and one has from equations (\ref{slowfast}):
\[
r\dot s_{1}+(1-r)\dot s_{2}+\epsilon\dot s_{{\rm{r}}}+\epsilon\dot
x_{{\rm{r}}}=-\epsilon \mu(s_{{\rm{r}}}^\star)x_{{\rm{r}}}<0 \ , 
\]
which implies that the trajectories are bounded.
For any $\sigma \geq \underline s$, $\hat s_{{\rm{r}}}^\star(\sigma)$ being the unique maximizer of the function $\beta(\sigma,\cdot)$, one has
\[
\sigma-\hat s_{{\rm{r}}}^\star(\sigma)=\frac{\mu(\hat s_{{\rm{r}}}^\star(\sigma))}{\mu^\prime(\hat s_{{\rm{r}}}^\star(\sigma))} \ .
\]
The function $\mu(\cdot)$ being increasing and concave, one obtains the
inequality
\[
\sigma-\hat s_{{\rm{r}}}^\star(\sigma)\geq \eta:=\frac{\mu(\hat s_{{\rm{r}}}^\star(\underline s))}{\mu^\prime(\hat s_{{\rm{r}}}^\star(\underline s))}>0 \, , \quad \forall \sigma \geq
\underline s \ .
\]
Furthermore, one can write
\[
r\dot s_{1} + \epsilon \dot
s_{{\rm{r}}}=-d\epsilon(s_{1}-s_{2})-\mu(s_{{\rm{r}}})x_{{\rm{r}}}<0 \ .
\]
Thus $r s_{1}+\epsilon s_{{\rm{r}}}$ is decreasing and has a limit when $t$
tends to $+\infty$.
Since the trajectories are bounded, $r\dot s_{1}+\epsilon\dot s_{{\rm{r}}}$ is
uniformly continuous, and we conclude by Barbalat's Lemma (see for instance \cite{K01}) that $r\dot s_{1} + \epsilon \dot s_{{\rm{r}}}$ converges
to $0$,
which implies that the positive quantities $s_{1}-s_{2}$ and $\mu(s_{{\rm{r}}})x_{{\rm{r}}}$ have to converge also to $0$. Notice that
$s_{{\rm{r}}}=0$ implies $\dot s_{{\rm{r}}}=\mu(s_{{\rm{r}}}^\star)s_{1}>\mu(\hat s_{{\rm{r}}}^\star(\underline s))\underline s>0$. So $s_{{\rm{r}}}$ cannot tend
to $0$ and $x_{{\rm{r}}}$ has necessarily to converge to $0$.
Write now the dynamics
\[
\frac{d}{dt}(s_{1}-s_{{\rm{r}}})=-\left(1+\frac{\epsilon}{r}\right)\mu(s_{{\rm{r}}}^\star)(s_{1}-s_{{\rm{r}}})-
d\frac{\epsilon}{r}(s_{1}-s_{2})-\mu(s_{{\rm{r}}})x_{{\rm{r}}} \ ,
\]
where $\mu(s_{{\rm{r}}}^\star)>\mu(\hat s_{{\rm{r}}}^{\star}(\underline s))>0$ and 
$d\frac{\epsilon}{r}(s_{1}-s_{2})-\mu(s_{{\rm{r}}})x_{{\rm{r}}}$
tends to $0$. Thus, there exists a time $T>0$ large enough such that 
\[
s_{{\rm{r}}}(t)>s_{1}(t)-\eta \geq s_{{\rm{r}}}^\star(t) \ , \quad  \forall t >T \ ,
\]
which implies to have $\mu(s_{{\rm{r}}})-\mu(s_{{\rm{r}}}^\star)>0$ for large $t$,
thus a contradiction with the convergence of $x_{{\rm{r}}}$ to $0$.\\

\medskip

Clearly the feedback (\ref{MRAP}) leaves the set $\{s_{1}=s_{2}\}$
invariant. Denote for simplicity $s_{l}=s_{1}=s_{2}$, and write
\[
\dot s_{l}+\epsilon\dot s_{{\rm{r}}}+\epsilon\dot x_{{\rm{r}}}=-\epsilon
\mu(s_{{\rm{r}}}^\star)x_{{\rm{r}}}<0 \ .
\]
Trajectories are thus bounded, and by Barbalat's Lemma one obtains
that $\mu(s_{{\rm{r}}}^\star)x_{{\rm{r}}}$ tends to $0$.
We prove now that $s_{l}$ has to reach $\underline s$ in
finite time. If not, $\mu(s_{{\rm{r}}}^\star(t))>\mu(\underline s)$ for any
time and $x_{{\rm{r}}}$
tends to zero. Write the dynamics
\[
\frac{d}{dt}(s_{l}-s_{{\rm{r}}}) =
-(1+\epsilon)\mu(s_{{\rm{r}}}^\star)(s_{l}-s_{{\rm{r}}})+\mu(s_{{\rm{r}}})x_{{\rm{r}}} \ .
\]
As before, we deduce that there exits a time $T^\prime>0$ such that
\[
s_{{\rm{r}}}(t)>s_{l}(t)-\eta\geq s_{{\rm{r}}}^\star(t) \, , \quad  \forall t> T^\prime \ ,
\]
leading to a contradiction with the convergence of $x_{{\rm{r}}}$ to $0$.
\qed

\end{document}